\documentclass{amsart}
\usepackage{amsfonts,amsmath,amsthm,amssymb}

\newtheorem{theorem}{Theorem}[section]

\newtheorem{lemma}{Lemma}[section]

\newtheorem{prop}{Proposition}[section]

\newtheorem{corollary}{Corollary}[section]

\renewcommand{\leq}{\leqslant}
\renewcommand{\geq}{\geqslant}

\begin{document}

\title{On 2-groups of Almost Maximal Class}
\author{Czes{\l}aw Bagi\'nski, Alexander Konovalov}
\date{}

\begin{abstract}
Let $G$ be a $2$-group of order $2^n$, $n\geq 6$, and nilpotency
class $n-2$. The invariants of such groups determined by their
group algebras over the field of two elements are given in the
paper.
\end{abstract}

\maketitle

\section{Introduction}
\label{Int}

Let $p$ be a prime. We say that a finite $p$-group $G$ of order
$p^n$ has the {\em coclass $k$} if the nilpotency class of $G$ is
equal to $n-k$. The $p$-groups of coclass $k=1$ and $k=2$ are
called the $p$-groups of {\em maximal class} and of {\em almost
maximal class} respectively. The natural partition of the class of
all $p$-groups into a series of well-behaved families of
$p$-groups of fixed coclass has proved to be one of the most
fruitful ideas of classifying $p$-groups (see
\cite{Leedham-GreenMcKay2002}). But the full description of
$p$-groups by generators and relations within each family  was
obtained only for a small number of families. It is well known for
$(p,k)\in\{(2,1),(3,1)\}$ and can be derived for $(p,k)\in\{
(2,2),(2,3)\}$ from \cite{Newman99} (see also \cite{James75},
\cite{James83}).

Here we study the internal structure of $2$-groups of nilpotency
coclass $2$ counting the important invariants of these groups. The
work, though motivated by the Modular Isomorphism Problem (MIP)
for finite $p$-groups is of independent interest. Recall that the
MIP asks whether, for given finite $p$-groups $G$ and $H$, the
isomorphism of the group algebras $FG$ and $FH$ over a field $F$
of characteristic $p$ implies the isomorphism of $G$ and $H$. The
problem was posed more than 50 years ago and for a long time there
was rather small progress in studying it. The state of the problem
in the end of the nineties in a more general context is briefly
described in \cite{Bleher} (see also \cite{Sand-survey} and
\cite{Baginski88a}, \cite{Baginski88b}, \cite{Salim-Sand},
\cite{Sand-metacyclic}, \cite{Wursthorn-diploma}, \cite{Wursthorn}).
One of the last results known to the
authors are that of Wursthorn saying that the groups of order
$2^6$ and $2^7$ are determined by their modular group algebras
over $GF(2)$ (see \cite{Bleher,Wursthorn}). He used computer
programs in studying the problem for groups of small order, which
are classified and have descriptions usable by computers. In December
2002 the authors were informed that the problem was solved
positively for group algebras of all $p$-groups over the field
$F=GF(p)$ (\cite{BorgeLaudal} \footnote{The proof given in the paper appeared
to be incorrect, so the MIP remains open.}).

Throughout, $p$ will denote a fixed prime. If $G$ is a $p$-group
then $G'=(G,G)$ is the commutator subgroup of $G,$ $\Phi(G)$ is
the Frattini subgroup, $G^p$ is the subgroup generated by $p$-th
powers of all elements of $G$ and $\Omega_i(G)$ is the subgroup
generated by all elements of order not bigger than $p^i$. By
$\gamma_i(G),$ $i\geq 2,$ we mean the $i$-th term of the lower
central series of $G$ ($\gamma_2(G)=G'$); by $\gamma_1(G)$ we mean
the subgroup $C_G(\gamma_2(G) \mod \gamma_4(G))$; if $g\in G$ and
$X\subset G$, then $(X,g)$ denotes the set of commutators
$\{(x,g):\ x\in X\}$. By $d(G)$ we denote the minimal number of
generators of $G,$ that is $d(G)=\log_p(|G/\Phi(G)|).$ The center
of $G$ we denote by $\zeta(G).$ If $g,h\in G$ then $g^h=h^{-1}gh$
and $C_g$ is the conjugacy class of $g$, i.e. $C_g=\{g^h:\ h\in G
\}$. The set of all conjugacy classes of $G$ will be denoted by
$Cl(G).$ If $S$ is a normal subset of $G$ then $Cl_G(S)$ is the
set of all conjugacy classes of $G$ contained in $S.$

We use $(x,y)= x^{-1}y^{-1}xy$ for the group commutator of
elements $x$ and $y$ of a group. The group commutator of weight
$n,$ $n>2,$ we define inductively by $(x_1,x_2,\ldots
,x_n)=((x_1,x_2,\ldots ,x_{n-1}),x_n).$

For a $p$-group $G$ we will denote by $R(G)$ the {\em Roggenkamp
parameter}
$$\sum_{i=1, \dots ,t} d(C_G(g_i)),$$ where
${g_1, \dots, g_t}$ is a set of representatives of all conjugacy
classes of $G$. This parameter was introduced by Roggenkamp and
described in \cite{Wursthorn-diploma}, \cite{Wursthorn}, where it
is proved that $R(G)$ is determined by the modular group algebra.
We will also use notation $R_G(S)$ for
$$\sum_{i=1, \dots ,s} d(C_G(h_i)),$$
where ${h_1, \dots, h_s}$ is a set of representatives of the conjugacy
classes of $G$ contained in the normal subset $S$ of $G$.

If $G$ is a finite $p$-group then by \cite{Quillen} the number of
conjugacy classes of maximal elementary abelian subgroups of given
rank is determined by $FG.$ By the {\em Quillen parameter} of $G$
we mean a series $Q(G)=(q_1,q_2,\ldots),$ where $q_i$ is the
number of conjugacy classes of maximal elementary abelian
subgroups of rank $i.$ Since $2$-groups of nilpotency coclass
$\leq 2$ do not contain elementary abelian subgroups of rank
bigger than $4,$ we will write the Quillen parameter in the form
$(q_1,q_2,q_3,q_4).$

In the paper we count the isomorphism class of the centers, the
numbers of conjugacy classes, the Quillen parameters and the
Roggenkamp parameters of all $2$-groups of nilpotency coclass $2$.
From this follows that the last two parameters determine almost
all the groups.

\begin{theorem}
Let $G$ and $H$ be groups of order $2^n,\ n\geqslant 8,$ and
nilpotency class $n-2.$ If $Q(G)=Q(H)$ and $R(G)=R(H)$ then either
$G\simeq H$ or
$\left\{G,H\right\}\subset\left\{G_9,G_{13},G_{14}\right\}$ or
$\left\{G,H\right\}=\left\{G_{24},G_{25}\right\}$.
\end{theorem}
It appears also that the numbers of conjugacy classes does not
give more information about differences between these groups and
the isomorphism class of the center allows only to exclude $G_9$
from the first set in the theorem.

Using the fact that every metacyclic $p$-group $G$ is determined
by the modular group algebra $FG$ over the field $F=GF(p)$
(\cite{Baginski88b},\cite{Sand-metacyclic}) and some additional
arguments one can prove the following theorem.

\begin{theorem}
Let  $F=GF(2)$ be the field of $2$ elements and let $G$ be a group
of order $2^n,\ n\geqslant 8,$ and nilpotency class $n-2.$ Then
$FG \simeq FH \Rightarrow G \simeq H$, provided
$\left\{G,H\right\}\neq\left\{G_{24},G_{25}\right\}$
\end{theorem}

As we were informed by W. Kimmerle for $n=8$ it was checked using
a computer that $FG_{24}$ is not isomorphic to $ FG_{25}.$

\section{Presentations of groups of almost maximal class}
\label{presentation}

Following \cite{James75}, \cite{James83} and \cite{Newman99},
2-groups of almost maximal class of order $2^{n}, n > 4,$ are
classified. However, the full list of presentations of these
groups by generators and relations was not published. In
\cite{James75} and \cite{James83} some groups were omitted. In
\cite{Newman99} there is given the list of pro-$2$ presentations
of pro-$2$ groups of coclass $\le 3$. Using both approaches first
we give the list of all groups of nilpotency coclass $2$ and order
$2^n,$ $n>4.$ According to \cite{Newman99} they are divided into
five families with numbers $7,8,9,50,59.$ We will use the same
numbers of the families and denote them by $Fam k,$ where
$k\in\{7,8,9,50,59\}.$ For our purposes the presentations derived
from the pro-$2$ presentations of pro-$2$ groups are more
convenient than ones that were given in James' paper.

\subsection{$2$-groups of almost maximal class with cyclic commutator subgroup.}

\begin{theorem}(\cite{James75},Theorem 5.1)\label{Fam59} There are precisely $6$ groups of
order $2^n$, $n\geqslant 4,$ and class $n-2$ with $G/\gamma_2(G)$
elementary abelian and $\gamma_2(G)$ cyclic. They form $Fam 59$ and are given by the
following presentation:\par
$$\langle x,y,t :\ x^{2^{n-2}}=t^2=1,\ y^2=z_1,\ x^y=x^{-1}z_2,\ x^t=xz_3,\  t^y=tz_4 \rangle $$
where the values of $z_i$, $1\leqslant i\leqslant 4,$ for
particular groups $G_m$, $1\leqslant m\leqslant 6$, are such as in
the following table.
$$\begin{array}{|l|c|c|c|c|c|c|}\hline
    & G_1               & G_2 & G_3       & G_4           & G_5         &      G_6     \\ \hline
z_1 & 1  & 1            & x^{2^{n-3}}     &  x^{2^{n-3}}  &  1          &  x^{2^{n-3}} \\
z_2 & 1  & x^{2^{n-3}}  & 1               & 1             &  1          & 1 \\
z_3 & 1  & 1            & 1               & 1             & x^{2^{n-3}} & x^{2^{n-3}} \\
z_4 & 1  & 1            & 1               & x^{2^{n-3}}   & 1           & 1 \\ \hline
\end{array}$$
\begin{center}
{\rm Table 1.}
\end{center}
\end{theorem}

\vspace{8pt}
\begin{theorem}(\cite{James75},Theorem 5.2)\label{Fam9} The number of groups
of order $2^n$, $n\geqslant 5,$ and class $n-2$ with $G/\gamma_2(G)$ of exponent $4$,
$\gamma_2(G)$ cyclic and $\gamma_1(G)/\gamma_2(G)$ cyclic is:
$$
\begin{cases}
3, & \text{if $n=5$,} \\
6, & \text{if $n>5$.}
\end{cases}
$$
They form $Fam 9$ and are given by the following presentation:
$$\langle x,y,t :\ x^{2^{n-2}}=t^2=1,\ y^2=z_1,\ x^y=x^{-1}z_2t,\ x^t=xz_3,\  t^y=tz_4 \rangle $$
where the values of $z_i$, $1\leqslant i\leqslant 4,$ for
particular groups $G_m$, $7\leqslant m\leqslant 12$, are such as
in the following table (for $G_9,\ G_{11}$ and $G_{12}$ we have
$n>5$).
$$\begin{array}{|l|c|c|c|c|c|c|}\hline
    & G_7  & G_8         & G_9         & G_{10}      & G_{11}      &   G_{12}     \\ \hline
z_1 & 1    & x^{2^{n-3}} & 1           & 1           & 1           &  x^{2^{n-3}} \\
z_2 & 1    & 1           & x^{2^{n-4}} & 1           & x^{2^{n-4}} & x^{2^{n-4}}  \\
z_3 & 1    & 1           & 1           & x^{2^{n-3}} & x^{2^{n-3}} & x^{2^{n-3}}  \\
z_4 & 1    & 1           & x^{2^{n-3}} & x^{2^{n-3}} & 1           & 1            \\ \hline
\end{array}.$$
\begin{center}
{\rm Table 2.}
\end{center}
\end{theorem}

\vspace{8pt}
\begin{theorem}(\cite{James75},Theorem 5.3(a))\label{Fam50} The number of groups of order
$2^n$, $n\geqslant 5$, and class $n-2$ with $G/\gamma_2(G)$ of
exponent $4,$ $\gamma_2(G)$ cyclic and $\gamma_1(G)/\gamma_2(G)$
elementary abelian is:
$$
\begin{cases}
3, & \text{if $n=5$,} \\
4, & \text{if $n>5$.}
\end{cases}
$$
They form $Fam 50$ and are given by the following presentation:
$$\langle x,y :\  x^{2^{n-2}}=1, \ y^4=z_1,  \   x^y=x^{-1}z_2 \rangle,$$
where the values of $z_1$ and $z_2$ are such as in the following table
(for $G_{16}$ we have $n>5$).
$$\begin{array}{|l|c|c|c|c|}\hline
    & G_{13}  & G_{14}        & G_{15}        & G_{16}      \\ \hline
z_1 &   1     &  1            &  x^{2^{n-3}}  &   1         \\
z_2 &   1     &  x^{2^{n-3}}  &  1            & x^{2^{n-4}} \\\hline
\end{array}.$$
\begin{center}
{\rm Table 3.}
\end{center}
\end{theorem}

\vspace{8pt}
\subsection{$2$-groups of almost maximal class with $2$-generated commutator subgroup.}

\begin{theorem}(\cite{James75},\cite{Newman99})\label{Fam8} The number of groups of order
$2^n,$ $n\geqslant 5,$ and class $n-2$ with $G/\gamma_2(G)$ of
exponent $4,$ $\gamma_2(G)$ $2$-generated and
$\gamma_1(G)/\gamma_2(G)$ elementary abelian is
$$
\begin{cases}
3, & \text{if $n=5$,} \\
4, & \text{if $n=6$,} \\
9, & \text{if $n>6$, $n$ odd,} \\
10,& \text{if $n>6$, $n$ even.}
\end{cases}
$$
They form $Fam 8.$
\end{theorem}

The description of the groups of $Fam 8$ by generators and
defining relations is more complicated than in the previous families.
Let $|G|=2^n=2^{2k+2+\epsilon}$,
$n\geqslant 6,$ $\epsilon\in\{0,1\}.$ Each of the groups of $Fam 8$ can be described as
the group
\begin{equation}\begin{array}{lllll}
\vspace{4pt}\langle x_1,x_2,y:& & & & \\
x_1^{2^{k+\epsilon}}=x_2^{2^{k}}=1, & y^4=t_1, & x_1^y=x_1x_2, & x_2^y=x_1^{-2}x_2^{-1}t_2, & x_2^{x_1}=x_2t_3\rangle,
\end{array}\end{equation}
where $t_1,t_3\in\{1,z_1\},$ $t_2\in\{1,z_1,z_2,z_1z_2\},$
$\epsilon\in\{0,1\}$ and:\par if $\epsilon =0,$ then
$z_1=x_2^{2^{k-1}},\ z_2=x_1^{2^{k-1}};$\par if $\epsilon =1,$
then $z_1=x_1^{2^{k}},\ z_2=x_2^{2^{k-1}}.$\par

\noindent The values of $t_i$ ($i=1,2,3$), for the groups $G_m,$
$18\leqslant m\leqslant 27,$ are given in the following table.
Note that for $\epsilon =1$ the groups $G_{24}$ and $G_{25}$ are
isomorphic.
$$
\begin{array}{|c|c|c|c|c|c|c|c|c|c|c|}
\hline
    & G_{18} & G_{19} & G_{20} & G_{21} & G_{22} & G_{23} & G_{24} & G_{25} & G_{26} & G_{27}  \\

\hline
t_1 &   1   &   z_1   &   z_1  &    1   &   z_1  &    1   &   z_1  &    1     &    1   &   z_1    \\
t_2 &   1   &    1    &   z_1  &    1   &    1   &   z_1  &   z_2  &  z_1z_2  &   z_2  & z_1z_2  \\
t_3 &   1   &    1    &    1   &   z_1  &   z_1  &    1   &    1   &    1     &   z_1  &   z_1    \\
\hline
\end{array}
$$
\begin{center}
Table 4.
\end{center}

For $n=5$ we have only the groups $G_{18},$ $G_{19}$ and $$G_{17}=
\langle x_1,x_2, y | x_1^8 = x_2^4 =1, y^4=x_1^4, x_1^y=y^2x_1
x_2, x_2^y=x_1^{-2}, x_2^{x_1}=x_2^{-1}x_1^{2}x_2 \rangle; $$ for
$n=6$ we have only the groups $G_{18},$ $G_{19},$ $G_{20}$ and
$G_{23}.$

Accordingly to the definitions from \cite{Newman99}, $G_{18}$ is
the mainline group of this family with the immediate descendants
$G_{19}-G_{22}$ (which are terminal) and $G_{23}$ (which have $4$
immediate terminal descendants $G_{24}-G_{27}$ if $\epsilon =0$,
and $3$ immediate terminal descendants $G_{24}-G_{27}$ with
$G_{24}\cong G_{25}$, if $\epsilon =1$).

\subsection{$2$-groups of almost maximal class with $3$-generated commutator subgroup.}

\begin{theorem}(\cite{James75},\cite{Newman99})\label{Fam7} The number of groups of order
$2^n,$ $n\ge 6,$ and class $n-2$ with $G/\gamma_2(G)$ of exponent
$4,$ $\gamma_2(G)$ $3$-generated and $\gamma_1(G)/\gamma_2(G)$
elementary abelian is
$$
\begin{cases}
2, & \text{if $n=6$,} \\
4, & \text{if $n>6$, $n$ odd,} \\
12,& \text{if $n>6$, $n$ even.}
\end{cases}
$$
They form $Fam 7.$
\end{theorem}

The groups of $Fam7$ have the most complicated description by generators
and defining relations among all $2$-groups of coclass 2.

Assume first that $|G|=2^n=2^{2k+2}$, $n\geqslant 6,$ that is
$G=G_m,$ where $28\leqslant m\leqslant 39.$ Each of these groups
can be described as the group
\begin{equation}
\begin{array}{l}
\vspace{4pt}\langle x_1, x_2, y: \\
x_1^{2^{k+1}}=1,\, x_2^{2^{k-1}}=x_1^{2^{k}},\, y^4=t_1,\, x_1^y=y^2x_1 x_2t_2,\, x_2^y=x_1^{-2}t_3,\, x_2^{x_1}=x_2^{-1}t_4 \rangle,
\end{array}
\end{equation}
where $t_1,t_2,t_3\in\{1,z_1\},$ $t_4\in\{1,z_1,z_2\},$
$z_1=x_1^{2^{k}}=x_2^{2^{k-1}},z_2=x_1^{2^{k-1}}x_2^{2^{k-2}}.$
The values of $t_i$ ($i=1,2,3,4$) for particular groups are given
in the following table: $$
\begin{array}{|c|c|c|c|c|c|c|c|c|c|c|c|c|}
\hline
    & G_{28} & G_{29} & G_{30} & G_{31} & G_{32} & G_{33} & G_{34} & G_{35} & G_{36} & G_{37} & G_{38} & G_{39} \\

\hline
t_1 &   1   &    1    &   z_1  &   z_1  &    1   &    1   &   z_1  &   z_1  &    1   &    1   &   z_1  & z_1 \\
t_2 &   1   &   z_1   &    1   &   z_1  &   z_1  &    1   &   z_1  &    1   &    1   &    1   &   z_1  &  1  \\
t_3 &   1   &    1    &    1   &    1   &    1   &    1   &    1   &    1   &    1   &   z_1  &    1   &  1  \\
t_4 &   1   &    1    &   z_1  &   z_1  &   z_1  &   z_1  &    1   &    1   &   z_2  &   z_2  &   z_2  & z_2 \\
\hline
\end{array}
$$
\begin{center}
Table 5.
\end{center}

Note that for $n=6$ we have only groups $G_{28}$ and $G_{29}$.

For further calculations we derive additional helpful relations.
For $28\leqslant m\leqslant 35$ in a group $G=G_m$ we have
\begin{equation}\label{eq1}
x_1^{y^{2}}=x_1^{-1}x_2t_1,\ (x_1^2)^{y^2}=x_1^{-2}t_4,\
(x_1^2)^{y}=x_2t_4,\ x_2^{y^2}=x_2^{-1}t_4.
\end{equation}
If $36\leqslant m\leqslant 39$ then
\begin{equation}\label{eq2}
x_1^{y^{2}}=x_1^{-1}x_2t_1t_3z_1,\ (x_1^2)^{y^2}=x_1^{-2}z_1z_2,\
(x_1^2)^{y}=x_2t_3z_2,\ x_2^{y^2}=x_2^{-1}z_2.
\end{equation}

\vspace{8pt} If $|G|=2^n=2^{2k+3}$, $n> 6,$ that is $G\in\{
G_{40}, G_{41}, G_{42}, G_{43}\},$ the groups of $Fam7$ can be
described in a little simpler way:
\begin{equation}\begin{array}{l}
\vspace{4pt}\langle x_1,x_2, y:\\
\ x_1^{2^{k+1}}=x_2^{2^k}=1,\ y^4=t_1,\ x_1^y= y^2 x_1 x_2,\ x_2^y=x_1^{-2},\ x_2^{x_1}=x_2^{-1}t_4  \rangle
\end{array}\end{equation}
where $t_1,t_4\in\{1,z\},\ z=x_1^{2^{k}}x_2^{2^{k-1}}.$ The values of
$t_1$ and $t_4$ for particular groups are given in the following
table:
$$
\begin{array}{|c|c|c|c|c|}
\hline
{}   &  G_{40} & G_{41} & G_{42} & G_{43} \\
\hline
t_1  &    1    &    z   &   1    &    z   \\
t_4  &    1    &    z   &   z    &    1   \\
\hline
\end{array}$$
\begin{center}
Table 6.
\end{center}
The relations (\ref{eq1}) are valid also for these groups.

The group $G_{40}$ is the mainline group with the immediate
descendant $G_{28}$ which is mainline and seven immediate
descendants $G_{29}$--$G_{35}$ which are terminal. The groups
$G_{41}$ and $G_{43}$ are terminal, $G_{42}$ is capable with four
immediate descendants $G_{36}$--$G_{39}$ which in turn are
terminal. All the groups $G_{40}$--$G_{43}$ are immediate
descendants of the mainline group $G_{28}$. \vspace{8pt}

\section{Groups with cyclic commutator subgroup}
\label{sect3}

In this section we describe properties of the groups defined
in Theorems \ref{Fam59}, \ref{Fam9} and \ref{Fam50}. All these
groups are extensions of a noncyclic group of class at most two
with a cyclic maximal subgroup by the group of order $2$.

Let $G=G_m$, $1\leqslant m\leqslant 16$. We let $A$ be the
subgroup of $G$ generated by the elements $x$ and $t$ (for
$m\in\{13,14,15,16\}$ we put $t=y^2$). The following lemma is an
easy observation obtained by a straightforward computation from
the presentations of groups. We assume that all the groups $G_m$
have order $p^n$, where $n\geqslant 5.$
\begin{lemma}\label{AinG1G12}
Let $G=G_m\in Fam 9\cup Fam50\cup Fam59$. Then:\vspace{5pt}\par
(a) $\gamma_2(G)=\begin{cases} \langle x^2\rangle & \text{if $G\in
Fam50\cup Fam59$}\\ \langle x^2t\rangle & \text{if $G\in Fam9$}
\end{cases}$\par
(b) If $i\geqslant 3$ then $\gamma_i(G)=\langle
x^{2^{i-1}}\rangle$;\vspace{5pt}\par (c) $\zeta(G_m)=\begin{cases}
\langle x^{2^{n-3}},t\rangle & \text{for $m\in\{1,2,3,7,8,13,14\}$}, \\
\langle x^{2^{n-4}}t \rangle & \text{for $m\in\{4,9\}$}, \\
\langle y^{2} \rangle & \text{for $m=15$}, \\
\langle x^{2^{n-3}} \rangle & \text{for
$m\in\{5,6,10,11,12,16\}$};
\end{cases}$\par
(d) $|G:A|=2$ and $\Omega_1(A)=\langle x^{2^{n-3}},t\rangle$ is an
elementary abelian normal subgroup of $G$ of order $4$;\par (e) If
$A$ is nonabelian then $\gamma_2(A)=\langle x^{2^{n-3}}\rangle,$
$\zeta(A)=\langle x^2\rangle,$ in particular the nilpotency class
of $A$ is not greater than 2.
\end{lemma}

For counting the Roggenkamp and Quillen parameters of our groups
we need information about the conjugacy classes of elements of
order $2$ lying outside $A.$

\begin{lemma}\label{G-A-Fam9-59}
Let $G=G_m\in Fam 9\cup Fam50\cup Fam59$.\par (a) The set
$G\setminus A$ splits into the following four conjugacy classes:
$C_y,$ $C_{yx},$ $C_{yt}$, $C_{yxt}.$\par (b) The orders of the
elements $y,yx,yt,yxt$ are given in the following tables.
$$\begin{array}{|l|cccccc||cccccc|}\hline
    & G_1 & G_2 & G_3 & G_4 & G_5 & G_6 & G_7 & G_8 & G_9 & G_{10} & G_{11} & G_{12} \\ \hline
y   &  2  &  2  &  4  &  4  &  2  &  4  &  2  &  4  &  2  &  2     &   2    &  4  \\ \hline
yx  &  2  &  4  &  4  &  4  &  2  &  4  &  4  &  4  &  8  &  4     &   8    &  8  \\ \hline
yt  &  2  &  2  &  4  &  2  &  2  &  4  &  2  &  4  &  4  &  4     &   2    &  4  \\ \hline
yxt &  2  &  4  &  4  &  2  &  4  &  2  &  4  &  4  &  8  &  4     &   8    &  8  \\ \hline
\end{array}$$
\begin{center}
{\rm Table 7.}
\end{center}
$$\begin{array}{|l|cccc|}\hline
             & G_{13}  & G_{14} & G_{15} & G_{16}  \\ \hline
y            &  4      &  4     &  8     &  4    \\ \hline
yx           &  4      &  4     &  8     &  8    \\ \hline
yt=y^{-1}    &  4      &  4     &  8     &  4    \\ \hline
yxt=y^{-1}x  &  4      &  4     &  8     &  8    \\ \hline
\end{array}$$
\begin{center}
{\rm Table 8.}
\end{center}
\end{lemma}
\begin{proof}
(a) The unique subgroup of order $4$ contained in $\langle
x\rangle$ is equal to $\langle x^{2^{n-4}}\rangle$. Since
$(x^{2^{n-4}})^y\neq x^{2^{n-4}}$ and
$(x^{2^{n-3}})^y=x^{2^{n-3}},$ we have $C_{\langle x\rangle}(y)=
\langle x^{2^{n-3}}\rangle$, that is $|C_{\langle
x\rangle}(y)|=2.$ But $\langle x^{2^{n-4}}\rangle\in\zeta(A)$, so
for all $g\in\{y,yx,yt,yxt\}$ we have $|C_{\langle
x\rangle}(g)|=2$ and then $|C_{\langle x\rangle}(g)|\leqslant 8,$
as $|G:\langle x\rangle |=4.$ Hence $|C_g|\geq 2^{n-3}$. On the
other hand the size of the conjugacy class $C_g$ does not exceed
the size of the commutator subgroup because $C_g\subseteq gG'$. In
our groups we have $|G'|=2^{n-3}$. Consequently $|C_g|=2^{n-3}$
 and $C_g=gG'$. Finally, since the elements $y,yx,yt,yxt$ lie
in different cosets of $G$ by $G'$, we obtain $G\setminus
A=yG'\cup yxG'\cup ytG'\cup yxtG'= C_y\cup C_{yx}\cup C_{yt}\cup
C_{yxt}.$\par
(b) By the defining relations (Theorem \ref{Fam59}) for the groups
$G_m\in Fam59$ we have $y^2=z_1$, $(yx)^2=y^2x^yx=z_1z_2,$
$(yt)^2=y^2t^yt=z_1z_4$ and $(yxt)^2=y^2x^yt^yxt=z_1z_2z_3z_4.$
Using the values of $z_i$, $1\leqslant i \leqslant 4,$ we obtain
the orders given in the table.
\par Similarly, for the groups
$G_m\in Fam9$ we have $y^2=z_1$, $(yx)^2=y^2x^yx=z_1z_2z_3t,$
$(yt)^2=y^2t^yt=z_1z_4$ and $(yxt)^2=y^2x^yt^yxt=z_1z_2z_4t$ and
again using the values of $z_i$, $1\leqslant i \leqslant 4,$ we
obtain the orders given in the table.\par Finally, for the groups
of $Fam50$ we have $y^4=(y^{-1})^4=z_1$ and
$(y^x)^4=z_1z_2^2=(y^{-1}x)^4.$ The orders of the considered
elements are easily seen from the defining relations of the groups
of $Fam50.$
\end{proof}

\begin{prop}\label{Fam9-59-Cl}
Let $G=G_m\in Fam 9\cup Fam50\cup Fam59$.\par (a) If $A$ is
abelian, that is $m\in\{1,2,3,4,7,8,9,13,14,15\},$ then
$$|Cl(G)|=2^{n-2}+6.$$\par (b) If $A$ is nonabelian, that is
$m\in\{5,6,10,11,12,16\},$ then
$$|Cl(G)|=5\cdot 2^{n-5}+6.$$
\end{prop}
\begin{proof}
(a) It follows immediately from the defining relations and Lemma
\ref{AinG1G12} that $A$ is abelian if and only if $|\zeta(G)|=4.$
Now, if $A$ is abelian and $a\in A\setminus \zeta(G)$, then
$C_G(a)=A$ that is $|C_a|=2$. Thus
$|Cl_G(A\setminus\zeta(G))|=\frac{|A|-|\zeta(G)|}{2}=\frac{2^{n-1}-2^2}{2}=2^{n-2}-2.$
Now the four conjugacy classes contained in $G\setminus A$, the
four $1$-element classes and the classes counted above give the
formula.\par (b) For nonabelian $A$ we have $\zeta(A)=\langle
x^2\rangle$. The elements of $\zeta(A)\setminus\zeta(G)$ have
centralizers equal to $A$. So $|Cl_G(\langle
x^2\rangle)|=\frac{|\zeta(A)|-|\zeta(G)|}{2}+2=\frac{2^{n-3}-2}{2}+2=2^{n-4}+1.$
The subgroup $\langle x^2,t\rangle$ has index $2$ in $A$. The set
$t\langle x^2\rangle=\langle x^2,t\rangle\setminus\langle
x^2\rangle$ contains exactly two elements of order $2$ (they form
a conjugacy class of the element $t$) and exactly two elements of
order $4$ (they form a conjugacy class of the element
$tx^{2^{n-4}}$). All other elements of this set form
$\frac{2^{n-3}-4}{4}=2^{n-5}-1$ four-element conjugacy classes.
Every conjugacy class contained in $A\setminus \langle
x^2,t\rangle$ has $4$ elements. Therefore we obtain
$|Cl(G)|=|Cl_G(G\setminus A)|+|Cl_G(A\setminus\langle
x^2,t\rangle)|+ |Cl_G(\langle x^2,t\rangle\setminus\langle
x^2\rangle)|+ |Cl_G(\langle
x^2\rangle)|+|\zeta(G)|=4+2^{n-4}+(2^{n-5}+1)+(2^{n-4}-1)+2=
2^{n-3}+2^{n-5}+6$\par
\end{proof}

\begin{prop}\label{R(G)G'-cyc}
Let $G=G_m\in Fam 9\cup Fam50\cup Fam59$. \par (a) If $A$ is
abelian, that is $m\in\{1,2,3,4,7,8,9,13,14,15\},$ then
$$R(G)=2^{n-1}+r_m,$$
where the values of $r_m$ are given in the following table.
$$\begin{array}{|l|cccc|ccc|ccc|}
\hline
     & G_1 & G_2 & G_3 & G_4 & G_7 & G_8 & G_9 & G_{13} & G_{14} & G_{15} \\ \hline
r_m  & 20  & 18  & 16  & 16  & 14  & 12  &  10  &  12     &  12  &  8  \\ \hline
\end{array}$$
\begin{center}
{\rm Table 9.}
\end{center}
(b) If $A$ is nonabelian, that is
$m\in\{5,6,10,11,12,16\},$ then
$$R(G)=2^{n-2}+r_m,$$
where the values of $r_m$ are given in the following table.
$$\begin{array}{|l|cc|ccc|c|}
\hline
     & G_5 & G_6 & G_{10} & G_{11} & G_{12} & G_{16} \\ \hline
r_m  & 18  & 16  &  13    &  13    &  11    &  10 \\ \hline
\end{array}$$
\begin{center}
{\rm Table 10.}
\end{center}
\end{prop}
\begin{proof}
(a) It is seen from Lemma \ref{AinG1G12} that $A$ is abelian if
and only if $|\zeta(G)|=4$. But $\zeta(G)\leqslant A,$ so for
$a\in A\setminus\zeta(G)$ we have $C_G(a)=A,$ which is
$2$-generated. For $a\in\zeta(G)$ we have obviously $C_G(a)=G$.
Therefore by the proof of Lemma \ref{Fam9-59-Cl}(a) we obtain
\begin{equation}\label{R(A)}
R_G(A)=2\cdot |Cl_G(A\setminus\zeta(G))|+4\cdot
d(G)=2^{n-1}+\begin{cases}
8 & \text{if $m\in\{1,2,3,4\}$}\\
4  & \text{if $m\in\{7,8,9,13,14,15\}$}
\end{cases}
\end{equation}
As it was already noted in the proof of Proposition
\ref{G-A-Fam9-59} for $g\in G\setminus A$ we have $|C_g|=2^{n-3}$,
so $|C_G(g)|=2^3$ and $C_G(g)=\langle g,\zeta(G)\rangle$. If $g$
has order $2$ and $\zeta(G)$ is elementary abelian, then
$d(C_G(g))=3$. If $\zeta(G)$ is cyclic of order $4$ (i.e.
$m\in\{4,9\}$) or $g$ has order $4$, then obviously $d(C_G(g))=2.$
Finally $d(C_G(g))=1$, when $g$ has order $8$. Now using the
formula (\ref{R(A)}) and the information from the tables of
Proposition \ref{G-A-Fam9-59} we get the assertion.\par (b) It is
obvious that for $a\in\zeta(A)\setminus\zeta(G)$ we have
$C_G(a)=A$ that is $d(C_G(a))=2$. Hence
$$R_G(\langle x^2\rangle)=2\cdot\frac{2^{n-3}-2}{2}+2\cdot d(G)=
2^{n-3}+\begin{cases} 4 & \text{if $m\in\{5,6\}$} \\
2 & \text{if $m\in\{10,11,12,16\}$}. \end{cases}$$ If the
conjugacy class of an element $a\in t\langle x^2\rangle$ has four
elements, then obviously its centralizer is equal to $\langle x^2,
a\rangle=\langle x^2, t\rangle,$ which is $2$-generated. For the
representatives $t$ and $x^{2^{n-4}}t$ of the $2$-element classes
(see the proof of Proposition \ref{Fam9-59-Cl}b) we have
$$C_G(t)=\begin{cases} \langle x^2,y,t \rangle & \text{if
$m\in\{5,6,11,12\}$,} \\ \langle x^2,xy \rangle & \text{if
$m\in\{10,16\}$}
\end{cases}$$
and
$$C_G(tx^{2^{n-4}})=\begin{cases} \langle x^2,y,t \rangle & \text{if
$m\in\{5,6,10\}$}, \\ \langle x^2,xy \rangle & \text{if
$m\in\{11,12,16\}$}.
\end{cases}$$
Thus $$R_G(t\langle
x^2\rangle)=2\cdot\frac{2^{n-3}-4}{4}+\begin{cases} 6 & \text{if
$m\in\{5,6\} $} \\ 5 & \text{if $m\in\{10,11,12\} $} \\ 4 &
\text{if $m=16$}.
\end{cases}$$
It is clear that all the elements $a\in A\setminus\langle
x^2,t\rangle$ have order $2^{n-2}$ and by Proposition
\ref{Fam9-59-Cl}b we have $|C_a|=4$. Thus $C_G(a)=\langle a
\rangle$ and, consequently, $R_G(A\setminus\langle
x^2,t\rangle)=|Cl_G(A\setminus\langle x^2,t\rangle)|=2^{n-4}.$
Finally,
$$R(G)=R_G(G\setminus A)+2^{n-2}+
\begin{cases} 8 & \text{if
$m\in\{5,6\} $}, \\ 5 & \text{if $m\in\{10,11,12\} $} \\ 4 &
\text{if $m=16$},
\end{cases}$$
so after counting the centralizers of the elements
$g\in\{y,yx,yt,yxt\}$ in a similar way as in the proof of part (a)
we get the assertion.
\end{proof}

We finish this section with counting the Quillen parameters. It
was noted in Lemma \ref{AinG1G12} that $\Omega_1(A)$ is an
elementary abelian normal subgroup of $G$ of order $4$. Hence, if
there are no elements of order $2$ outside $A$, $\Omega_1(A)$ is
the unique maximal elementary abelian subgroup of $G$, i.e. the
Quillen parameter of $G$ is equal to $(0,1,0,0)$. Since
$|G:A|=2$, there is not an elementary abelian subgroup of order
$4$ with trivial intersection with $A$. Hence, in all other cases
the parameter has type $(0,q_2,q_3,0)$. The  values of $q_2$ and
$q_3$ can be easily counted by studying the centralizers of the
elements of order $2$ lying outside $A$ and the intersections of
these centralizers with $\Omega_1(A)$. We leave the details to the
reader.

\begin{prop}
If $G\in Fam9\cup Fam50\cup Fam59$ then the Quillen parameters
are given in the following table.$$
\begin{tabular}{|c|c|c|c|}
\hline
{\rm Group}    & {\rm Isomorphism}          & {\rm Quillen}      & {\rm Representatives of the conjugacy classes}   \\
               & {\rm type of $\zeta(G)$}   & {\rm parameter}    & {\rm of maximal elementary abelian subgroups} \\ \hline
$G_1$    & $C_2 \times C_2$     & $(0,0,2,0)$  & $ \langle
x^{2^{n-3}}, y, t \rangle, $\ $\langle x^{2^{n-3}}, yx, t \rangle
$ \\ \hline $G_2$    & $C_2 \times C_2$     & $(0,0,1,0)$  & $
\langle x^{2^{n-3}}, y, t \rangle $ \\ \hline $G_3$    & $C_2
\times C_2$     & $(0,1,0,0)$  & $ \langle x^{2^{n-3}},    t
\rangle $ \\ \hline $G_4$    & $C_4$                & $(0,3,0,0)$
& $ \langle x^{2^{n-3}}, y    \rangle, $\ $\langle x^{2^{n-3}}, yx
\rangle, $ $\langle x^{2^{n-3}}, x^{2^{n-4}}t \rangle $ \\ \hline
$G_5$    & $C_2$                & $(0,1,1,0)$  & $\langle
x^{2^{n-3}}, yx \rangle,$ $ \langle x^{2^{n-3}}, y, t \rangle$ \\
\hline
$G_6$    & $C_2$                & $(0,2,0,0)$  & $ \langle x^{2^{n-3}},    t \rangle $, $ \langle x^{2^{n-3}},    xyt \rangle $ \\
\hline\hline
$G_{7} $ & $C_2 \times C_2$     & $(0,0,1,0)$  & $ \langle x^{2^{n-3}}, y, t \rangle $ \\ \hline
$G_{8} $ & $C_2 \times C_2$     & $(0,1,0,0)$  & $ \langle x^{2^{n-3}},    t \rangle $ \\ \hline
$G_{9} $ & $C_4$                & $(0,2,0,0)$  & $ \langle x^{2^{n-3}},    t \rangle,$\ $\langle x^{2^{n-3}} , y \rangle $ \\ \hline
$G_{10}$ & $C_2$                & $(0,2,0,0)$  & $ \langle x^{2^{n-3}},    t \rangle,$\ $ \langle x^{2^{n-3}} , y \rangle $ \\ \hline
$G_{11}$ & $C_2$                & $(0,0,1,0)$  & $ \langle x^{2^{n-3}}, y, t \rangle $ \\ \hline
$G_{12}$ & $C_2$                & $(0,1,0,0)$  & $ \langle x^{2^{n-3}},    t \rangle $ \\ \hline\hline
$G_{13}$ & $C_2 \times C_2$     & $(0,1,0,0)$  & $ \langle x^{2^{n-3}},    y^2 \rangle,$ \\ \hline
$G_{14}$ & $C_2 \times C_2$     & $(0,1,0,0)$  & $ \langle x^{2^{n-3}},    y^2 \rangle,$ \\ \hline
$G_{15}$ & $C_4$                & $(0,1,0,0)$  & $ \langle x^{2^{n-3}}, y^2x^{2^{n-4}}\rangle$ \\ \hline
$G_{16}$ & $C_2$                & $(0,1,0,0)$  & $ \langle x^{2^{n-3}},    y^2 \rangle $ \\ \hline\hline
\end{tabular}
$$
\begin{center}
{\rm Table 11.}
\end{center}
\end{prop}

\section{Groups with 2-generated commutator subgroup} \label{family8}

In this section we describe the groups defined in Theorem
\ref{Fam8}. All these groups are extensions of a subgroup of
nilpotency class $\leqslant 2$ by the cyclic group of order $4$.

Let $G=G_m$, $18\leqslant m\leqslant 27$. In this section we let
$A$ be the subgroup of $G$ generated by the elements $x_1$ and
$x_2$. The following lemmas are easy observations obtained by a
straightforward computation from the presentations of groups. We
assume that all the groups $G_m$ have order $p^n$, where
$n\geqslant 7.$

\begin{lemma} \label{fam8-lcs}
Let $G \in Fam8$. Then $\gamma_2(G)= \langle x_1^2, x_2   \rangle
$, $\gamma_3(G)= \langle x_1^2, x_2^2 \rangle $, and in general,
$\gamma_{2i}(G)   = \langle x_1^{2^i}, x_2^{2^{i-1}} \rangle $,
$\gamma_{2i+1}(G) = \langle x_1^{2^i}, x_2^{2^i}     \rangle $ for
$i\geqslant 1$.
\end{lemma}

\begin{lemma}\label{fam8-A}
Let $ G \in Fam8$. Then:\par (a) If $G=G_m,$ where
$m \in \{ 21, 22, 26, 27 \},
$ then $A'=\{z_1\}$ and $\zeta(A)=\gamma_3(G);$ in
all other cases $A$ is abelian.\par (b) If
$g\in\{y,yx_1,y^{-1},y^{-1} x_1, y^2 x_1\}$, then
$(G,g)=\gamma_2(G);$ if $g\in\{y^2,y^2x_2\}$, then
$(G,g)=\gamma_3(G).$\par (c) $\Omega_1(A)$ is an elementary
abelian subgroup of order $4$ and $\Omega_1(A)\leq\zeta(\langle
A,y^2\rangle).$
\end{lemma}

\begin{lemma}\label{G-Acc}
The set $G \setminus A $ splits into 7 conjugacy classes, which
are
$ y          \gamma_2(G)$,
$ y      x_1 \gamma_2(G)$,
$ y^{-1}     \gamma_2(G)$,
$ y^{-1} x_1 \gamma_2(G)$,
$ y^2    x_1 \gamma_2(G)$,
$ y^2        \gamma_3(G)$,
$ y^2    x_2 \gamma_3(G)$.
Moreover, first four classes do not contain elements of order 2.
\end{lemma}

\begin{proof}
It is obvious that $A$ is a normal subgroup of index $4$ with the
cyclic factor group $G/A=\langle yA\rangle.$ We have also by Lemma
\ref{fam8-lcs} that $A=x_1\gamma_2(G)\cup
x_2\gamma_3(G)\cup\gamma_3(G).$ Hence
$$\begin{array}{ll}
G\setminus A= & yA\cup y^{-1}A\cup y^2A= \\
& y          \gamma_2(G)\cup
  y      x_1 \gamma_2(G)\cup
  y^{-1}     \gamma_2(G)\cup
  y^{-1} x_1 \gamma_2(G)\cup \\
& y^2    x_1 \gamma_2(G)\cup
  y^2        \gamma_3(G)\cup
  y^2    x_2 \gamma_3(G).
\end{array}$$
By Lemma \ref{fam8-A}(b) each of the sets from the last two lines
is a conjugacy class. Since $yA$ and $y^{-1}A$ are elements of order
$4$ in $G/A$, none of these cosets contains elements of order $2$
and the lemma is proved.
\end{proof}

\begin{lemma} \label{fam8-max-el-ab-sbgrp}
Let $ G \in Fam8$, and let $C$ be one of the conjugacy classes
$y^2 x_1 \gamma_2(G)$, $ y^2\gamma_3(G)$, $ y^2 x_2 \gamma_3(G).$
If $C$ consists of elements of order 2, then the family $\Theta =
\{ \langle g, x_1^{2^{k+\epsilon-1}}, x_2^{2^{k-1}} \rangle: g \in
C \}$ is a conjugacy class of maximal elementary abelian subgroups
of $G$. Moreover $|\Theta|=|C|/4$.
\end{lemma}

\begin{proof}
Notice first that any two elements from different conjugacy
classes listed in the lemma are not commuting. So they cannot lie
together in an abelian subgroup. Further, if $C$ is one of these
classes and $g\in C$, then as it follows from Lemma
\ref{fam8-A}(c), $\langle g,\Omega_1(G)\rangle$ is abelian. So by
the above, if $g$ has order $2$, it is maximal elementary abelian.
It is clear that one such subgroup contains $4$ elements belonging to $C$
and different subgroups determine disjoint such four-element
subsets. Thus $|\Theta|=\frac{|C|}{4}.$
\end{proof}

\begin{corollary} \label{fam8-order-max-el-ab-sbgrp}
If $ G \in Fam8$ and $|G| \ge 2^6$, then maximal elementary
abelian subgroups of $G$ have order not greater than $2^3$. If
$|G| > 2^6$, and $B$ is a maximal elementary abelian subgroup of
$G$ of order $2^3$, then $B$ is not normal in $G$.
\end{corollary}

The proof of the following lemma needs standard and not difficult
calculations.

\begin{lemma} \label{fam8-centralizers}
Let $ G \in Fam8$. Then: \par
(a) for $g \in \zeta(G)$, $C_G(g)=G$; \par
(b) for $g \in \langle z_1, z_2 \rangle \setminus \zeta(G)$,
    $C_G(g)= \langle y^2, x_1, x_2 \rangle $; \par
(c) for $g \in \gamma_3(G) \setminus \langle z_1, z_2 \rangle$,
    $C_G(g)= A $;\vspace{5pt} \par
(d) for $g \in A \setminus \gamma_3(G) $, $C_G(g)=\begin{cases} A
, & \text{ if $A$ is abelian,} \\ \langle g, \gamma_3(G) \rangle,
& \text{ if $A$ is nonabelian;} \end{cases}$ \vspace{5pt}\par (e)
for $g \in y^2 x_1 \gamma_2(G)$, $C_G(g)= \langle g, z_1, z_2
\rangle$; \par (f) for $g \in y^2 \gamma_3(G) \cup y^2 x_2
\gamma_3(G) $,
    $C_G(g)= \langle h, z_1, z_2 \rangle$, where $h^2=g$;\par
(g) for $g\notin \langle y^2,A\rangle,$ $C_G(g)=\langle
g,\zeta(G)\rangle.$
\end{lemma}

\begin{prop} \label{fam8-number-con-classes}
Let $ G \in Fam8$. \par (a) If $A$ is abelian, that is
$m\in\{18,19,20,23,24,25\},$ then $$|Cl(G)| = 9 +
2^{2k+\epsilon-2}.$$
\par (b) If $A$ is nonabelian, that is $m\in\{21,22,26,27\},$
then $$|Cl(G)| = 9 + 5 \cdot
2^{2k+\epsilon-5}.$$
\end{prop}

\begin{proof}
Assume first that $A$ is abelian. We split $A$ into the set
theoretic sum of $3$ disjoint subsets: $A=\zeta(G)\cup(\langle
z_1,z_2\rangle\setminus\zeta(G))\cup(A\setminus \langle
z_1,z_2\rangle).$ In $\zeta(G)$ we have two one-element classes,
the set $\langle z_1, z_2 \rangle \setminus \zeta(G)$ forms one
two-element class and finally in the last subset we have
$\frac{|A|-4}{4}=2^{2k -2+\epsilon}-1$ four-element classes. All
this classes together with the seven classes contained in
$G \setminus A$ give $9 + 2^{2k+\epsilon-2}$ conjugacy classes.

Now let $A$ be nonabelian and let us split it into the set theoretic
sum of $4$ disjoint sets $A=\zeta(G)\cup(\langle
z_1,z_2\rangle\setminus\zeta(G))\cup(\gamma_3(G)\setminus \langle
z_1,z_2\rangle)\cup(A\setminus\gamma_3(G)).$ As in the previous
case, conjugacy classes contained in the three first sets have
respectively $1,$ $2$ and $4$ elements. Each conjugacy class
contained in $A\setminus\gamma_3(G)$ has $8$ elements, by Lemma
\ref{fam8-centralizers}(d). Therefore
$|Cl(G)|=2+1+\frac{2^{2k+\epsilon-2}-4}{4}+\frac{2^{2k+\epsilon}-2^{2k+
\epsilon-2}}{8}+7= 9+5\cdot 2^{2k+\epsilon-5}.$
\end{proof}

\begin{lemma}
Let $G\in Fam8.$ The conjugacy classes of $G$ not contained in
$A$ have elements of order as listed in the following table.
$$
\begin{tabular}{|c|c|c|c|c|c|c|c|c|c|c|}
\hline {\rm Groups}         & $G_{18}$ & $G_{19}$ & $G_{20}$ & $G_{21}$ & $G_{22}$ & $G_{23}$  & $G_{24}$ & $G_{25}$  & $G_{26}$ & $G_{27}$ \\ \hline
$A$ {\rm abelian}           &    $+$   &   $+$    &   $+$    &   $-$    &   $-$    &   $+$     &   $+$    &   $+$     &    $-$   &    $-$   \\ \hline
{\rm Representatives of}    &          &          &          &          &          &           &          &           &          &          \\
{\rm conjugacy classes}     &          &          &          &          &          &           &          &           &          &          \\
$y,$ $y^{-1}$               &   $4$    &   $8$    &  $8$     &  $4$     &    $8$   &   $4$     &  $8$     &  $4$      &  $4$     &  $8$     \\
$yx_1,$ $y^{-1}x_1$         &   $4$    &   $8$    &  $8$     &  $8$     &    $4$   &   $4$     &  $4$     &  $8$      &  $4$     &  $8$     \\
$y^2x_1$                    &   $2$    &   $4$    &  $4$     &  $2$     &    $4$   &   $2$     &  $4$     &  $2$      &  $2$     &  $4$     \\
$y^2   $                    &   $2$    &   $4$    &  $4$     &  $4$     &    $2$   &   $2$     &  $2$     &  $4$      &  $2$     &  $4$     \\
$y^2x_2$                    &   $2$    &   $4$    &  $2$     &  $2$     &    $4$   &   $4$     &  $4$     &  $4$      &  $4$     &  $4$     \\
\hline
\end{tabular}
$$
\begin{center}
{\rm Table 12.}
\end{center}
\end{lemma}

\begin{proof}
The order of $y$, $y^{-1}$ and $y^2$ can be easily fixed by
Table 4. For $yx_1$ and $y^{-1}x_1$ we have
$(yx_1)^2=y^2x_1^2x_2^{x_1}\in y^2x_2\gamma_3(G)$ and similarly
$(y^{-1}x_1)^2\in y^2x_2\gamma_3(G).$ Further
$(y^2x_2)^2=t_1t_2t_2^yt_3.$ If $t_2=z_1$, then obviously
$t_2^y=t_2.$ If $t_2=z_2$, then $t_2^y=z_1z_2.$ Now it is an easy
task to fill the table.
\end{proof}

\begin{lemma} \label{fam8-Rog-of-A-in-G}
Let $G \in Fam8$.\par (a) If $A$ is abelian, that is
$m\in\{18,19,20,23,24,25\},$ then
$$R_G(A)=5+2^{2k+\epsilon-1}.$$ \par (b) If $A$ is nonabelian, that
is $m\in\{21,22,26,27\},$ then
$$R_G(A)=5+5 \cdot 2^{2k+\epsilon-4}.$$
\end{lemma}

\begin{proof}
The only conjugacy class contained in $A$ whose representatives have
$3$-generated centralizers is $C_{z_2}=\{z_2,z_1z_2\}.$ Representatives of all other
conjugacy classes contained in $A$ have $2$-generated
centralizers.
\end{proof}

Now we are ready to count the Roggenkamp
and the Quillen parameters of all the groups of $Fam8.$

\begin{prop}\label{Fam8-R}
Let $G=G_m\in Fam8,$ $|G|>2^6$. Then
$$R(G)=\begin{cases} 2^{2k+\epsilon-1}+r_m, & \text{if $A$ is abelian, that is $m\in\{18,19,20,23,24,25\},$} \\
               5\cdot2^{2k+\epsilon-4}+r_m, & \text{if $A$ is nonabelian, that is $m\in\{21,22,26,27\},$}
\end{cases}$$
where the values of $r_m$ are given in the following table.
$$\begin{array}{|l|cccccc|cccc|}
\hline
     & G_{18} & G_{19} & G_{20} & G_{23} & G_{24} & G_{25}  & G_{21} & G_{22} & G_{26} & G_{27} \\ \hline
r_m  & 20     & 15     & 16     & 19     & 17     & 17      & 18     & 17     &  19    &  15 \\ \hline
\end{array}$$
\begin{center}
{\rm Table 13.}
\end{center}

\end{prop}
\begin{proof}
By Lemma \ref{fam8-Rog-of-A-in-G} in order to find the Roggenkamp
parameters we need to count minimal numbers of generators of the
centralizers of the elements listed in Table 12. The elements
$y,y^{-1},yx_1,y^{-1}x_1$
have centralizers of order $8$, so if any of these elements has
order $4$, then it must have a $2$-generated centralizer. If it has
order $8$, then its centralizer is cyclic. The elements $y^2,
y^2x_2$ are $2$-powers and have centralizers of order $16$ which
are $2$-generated. The centralizer of $y^2x_1$ has $8$ elements
and depending on its order it is $3$-generated when its order is
equal $2$ and $2$-generated, when its order is equal $4$.
\end{proof}

\begin{prop}\label{Fam8-Q}
If $G=G_m\in Fam8,$ $|G|>2^6$, then the Quillen parameter of $G$
is such as it is listed in the following table.
$$
\begin{tabular}{|c|c|c|}
\hline
            & {\rm Quillen}   &  {\rm The representatives of the conjugacy classes} \\
{\rm Group} & {\rm parameter} &  {\rm of elementary abelian maximal subgroup}  \\ \hline
$G_{18}$ & $(0,0,3,0)$ & $\langle y^2x_1,\Omega_1(A) \rangle,\ \langle y^2,\Omega_1(A) \rangle,\ \langle y^2x_2,\Omega_1(A) \rangle$ \\ \hline
$G_{19}$ & $(0,1,0,0)$ & $\Omega_1(A)$ \\ \hline
$G_{20}$ & $(0,0,1,0)$ & $\langle y^2x_2,\Omega_1(A) \rangle$ \\ \hline
$G_{21}$ & $(0,0,2,0)$ & $\langle y^2x_1,\Omega_1(A) \rangle,\ \langle y^2x_2,\Omega_1(A) \rangle$ \\ \hline
$G_{22}$ & $(0,0,1,0)$ & $\langle y^2,\Omega_1(A) \rangle$ \\ \hline
$G_{23}$ & $(0,0,2,0)$ & $\langle y^2x_1,\Omega_1(A) \rangle,\ \langle y^2,\Omega_1(A) \rangle$ \\ \hline
$G_{24}$ & $(0,0,1,0)$ & $\langle y^2,\Omega_1(A) \rangle$ \\ \hline
$G_{25}$ & $(0,0,1,0)$ & $\langle y^2x_1,\Omega_1(A) \rangle$ \\ \hline
$G_{26}$ & $(0,0,2,0)$ & $\langle y^2x_1,\Omega_1(A) \rangle,\ \langle y^2,\Omega_1(A) \rangle$ \\ \hline
$G_{27}$ & $(0,1,0,0)$ & $\Omega_1(A)$ \\
\hline
\end{tabular}
$$
\begin{center}
{\rm Table 14.}
\end{center}
\end{prop}

\begin{proof}
It was mentioned earlier that if the elements $y^2,y^2x_1,y^2x_2$
have order $4$ then $\Omega_1(A)$ is the unique maximal elementary
abelian subgroup of $G$ and then $(0,1,0,0)$ is the Quillen
parameter of $G$. If among these elements there are $i$ elements
of order $2,$ $i\in\{1,2,3\},$ then the Quillen parameter of $G$
has the form $(0,0,i,0).$ So using the information from Table 12
one can easily fill Table 13.
\end{proof}

\section{Groups with 3-generated commutator subgroup} \label{family7}

In this section we describe the groups defined in Theorem
\ref{Fam7}.

Let $G=G_m$, $28\leqslant m\leqslant 43$. In this section we let
$A$ be the subgroup of $G$ generated by the elements $x_1^2$ and
$x_2$. This subgroup will play a key role in our computations.

The following lemmas are easy observations obtained by a
straightforward computation from the presentations of the groups. We
assume that all the groups $G_m$ have order $p^n$, where
$n\geqslant 6.$

\begin{lemma} \label{fam7-lcs}
Let $ G \in Fam7$. Then $\gamma_2(G)= \langle y^2x_1^2,
x_1^2x_2, x_2^2\rangle $, $\gamma_3(G)= \langle x_1^2x_2, x_2^2
\rangle $, and, in general for $m\geq 2,$ $\gamma_{2m}(G) = \langle x_1^{2^m},
x_2^{2^{m-1}} \rangle $, $\gamma_{2m+1}(G) = \langle x_1^{2^m}x_2^{2^{m-1}},
x_2^{2^m}\rangle $.
\end{lemma}

\begin{lemma}\label{fam7-A}
Let $G \in Fam7$. Then: \par (a) $\gamma_3(G)\le A$ and
$|A:\gamma_3(G)|=2;$ \par (b) $A$ is normal in $G$ and the factor group
$\overline{G}=G/A$ is isomorphic to the dihedral group of order
$8;$\par (c) If $G=G_m,$ with $36\le m\le 39,$ then $[A,A]=\langle
z_1\rangle;$ in all other cases $A$ is abelian;\par (d) The
subgroup $H=\langle y^2, A\rangle$ is the unique normal subgroup
of $G$ of index $4$ containing $A;$\par (e) There exist exactly
$4$ non-normal subgroups of $G$ of index $4$ containing $A$:
$H_1=\langle x_1, A\rangle ,$ $H_2=\langle y^2x_1, A\rangle ,$
$H_3=\langle yx_1^{-1}, A\rangle ,$ $H_4=\langle y^3x_1, A\rangle
.$ Moreover, $H_1^y=H_2$ and $H_3^{x_1}=H_4.$ \par
\end{lemma}

Let $H=\Phi(G)=\langle y^2, A\rangle.$ Since $d(G)=2,$ $G$ has
exactly three maximal subgroups. These are $M_1=\langle
y,H\rangle,$ $M_2=\langle x_1,H\rangle $ and $M_3=\langle
yx_1,H\rangle .$  It can be easily seen that $d(M_1)=d(M_2)=2$
and $d(M_3)=3.$ In further considerations we will use the
splitting of $G$ into the following set-theoretic sum of pairwise
disjoint subsets: $$G=A\cup (H\setminus A)\cup (M_1\setminus
H)\cup (M_2\setminus H)\cup (M_3\setminus H).$$

\begin{lemma} \label{fam7-M1}
Let $G \in Fam7.$ Then \par
(a) $Cl_G(H\setminus A)=\{ y^2\gamma_3(G), y^2x_1^{-2}\gamma_3(G)\},$ \par
(b) $Cl_G(M_1\setminus H)=\{ y\gamma_2(G),$ $y^3\gamma_2(G)\}.$ \par
\end{lemma}

\begin{proof}
(a) Since $H\setminus A=y^2A$ and $A=\gamma_3(G)\cup
x_1^{-2}\gamma_3(G)$, we have $H\setminus A=y^2\gamma_3(G)\cup
y^2x_1^{-2}\gamma_3(G)$. Now the relations \ref{eq1} and \ref{eq2} give $[G,y^2]=
[G,y^2x_1^{-2}]=\gamma_3(G).$ Hence for $g\in \{y^2, y^2x_1^{-2}\},$
$C_g=g\gamma_3(G).$


(b) By Lemma \ref{fam7-lcs}, $H=\gamma_2(G)\cup x_1^2\gamma_2(G),$
so $M_1\setminus H=yH=y\gamma_2(G)\cup yx_1^2\gamma_2(G).$ Now
$C_G(y)=\langle y,\zeta(G)\rangle$ and
$C_G(yx_1^2)=\langle yx_1^2, \zeta(G) \rangle$,
that is for $g\in \{y, yx_1^2 \}$,
$|C_x|=|\gamma_2(G)|$ and the assertion follows.\end{proof}

The following lemma will allow us to count the Quillen parameters of all
groups of $Fam7$.

\begin{lemma}\label{fam7-odd-order-2}
Let $G=G_m$ be a group of $Fam7.$ Then the set of all elements of order $2$ of $G$
is contained in
$$\Omega_1(A)\cup\bigcup_{g\in X}C_g$$
where
$$X=\left\{
\begin{array}{lll}
\{y^2,y^2x_1^{-2},yx_1^{-1}, yx_1^{-1}x_2^{2^{k-1}}\} & {\it if \ {\it A}\ is \ abelian,}\\
\{y^2, yx_1^{-1}x_2^{2^{k-2}}\} & {\it if \ {\it A} \ is \ nonabelian.}
\end{array}\right.$$
\end{lemma}

\begin{proof}
It suffices to show, that conjugacy classes which are not
represented by elements from $X,$ do not consist of elements of
order $2.$ By Lemma \ref{fam7-M1} elements belonging to
$M_1\setminus H$ are conjugated either to $y$ or to $y^3,$ so
their order is equal either $4$ or $8.$ It can be also easily
verified that elements of the set $M_2\setminus H$ have order
equal to $o(x_1)=2^{k+1}$. So let us consider elements of
$M_3\setminus H.$ Since $M_3\setminus H=yx_1^{-1}A\cup y^3x_1A$
and $(yx_1^{-1}A)^y=y^3x_1A$, we consider only elements from
$yx_1^{-1}A.$ Let $g=yx_1^{-1}(x_1^{2r}x_2^s)$ be an arbitrary
element of this set. First let us assume that $|G|=2^{2k+2}$ (in
this case $0\le r<2^k,\ 0\le s<2^{k-1}$). If $A$ is abelian, then
$g^2=t_1t_2t_4(x_1^{-2}x_2)^{s-r}.$ Since $t_1,t_2,t_4\in
\{1,x_1^{2^{k}}=x_2^{2^{k-1}}\},$ $g$ has order $2$ if and only if
$t_1t_2t_4=1=(x_1^{-2}x_2)^{s-r}.$ It follows from Table 5 that
for $m\in \{29,31,33,35\}$ $t_1t_2t_4=z_1$ so the subset
$M_3\setminus H$ of $G_m$ does not contain elements of order $2.$
If $m\in\{28,30,32,34\}$ then $t_1t_2t_4=1$ and
then $g$ is of order $2,$ when $r-s\equiv 0(\rm mod\ 2^{k-1}).$
Since this congruence has exactly $2^k$ solutions, there exist
exactly $2^{k+1}$ elements of order $2$ in $M_3\setminus H$ (half
of them lie in $yx_1^{-1}A$ and the second half in $y^3x_1A$). It
can be easily checked that for $m\in \{28,30\}$,
$C_G(yx_1^{-1})=\langle yx_1^{-1}, y^2x_1^{-2}, z_1,
x_1^{-1}x_2\rangle$, so $|C_{yx_1^{-1}}|=2^k.$ The second
conjugacy class consisting of elements of order $2$ is the class
represented by $yx_1^{-1}z_1.$ If $m\in\{32,34\}$, then
$C_G(yx_1^{-1})=\langle yx_1^{-1}, z_1, x_1^{-1}x_2 \rangle$ and
then all elements of order $2$ of $M_3\setminus H$ lie in one
conjugacy class.

Now let us assume that $A$ is nonabelian, that is
$m\in\{36,37,38,39\}.$ For $g=yx_1^{-1}(x_1^{2r}x_2^s)$ we have
\begin{equation}\label{nab-g2}
g^2=t_1t_2t_3^{1+r+s}z_1^{1+r(s-r)}z_2x_1^{-2(s-r)}x_2^{s-r}.
\end{equation}
Therefore, if $g^2=1$ we have $x_1^{-2(s-r)}x_2^{s-r}\in
\Omega_1(A),$ which says that
$x_1^{-2(s-r)}x_2^{s-r}\in\{1,x_1^{-2^{k-1}}x_2^{2^{k-2}}=z_1z_2\}.$
In particular $s$ and $r$ have to be of the same parity, and
because of that $g^2=t_1t_2t_3z_1z_2x_1^{-2(s-r)}x_2^{s-r}.$ Thus,
as it follows from Table 5, there is no element of order $2$ in
$M_3\setminus H$ for $G_{37}$ and $G_{39}.$ If $G=G_{36}$ or
$G=G_{38}$ then $g^2=1$ if and only if $r-s\equiv 2^{k-2}({\rm
mod}\ 2^{k}).$ This congruence has $2^k$ solutions which means, in
particular, that in $M_3\setminus H$ there exist exactly $2^{k+1}$
elements of order $2$ and they all belong to the class of the
element $yx^{-1}x_2^{2^{k-2}}$ as
$C_G(yx^{-1}x_2^{2^{k-2}})=\langle yx^{-1}x_2^{2^{k-2}}, z_1,
x_1^{-2}x_2\rangle$, i.e. $|C_{yx^{-1}x_2^{2^{k-2}}}|=2^{k+1}.$

Finally, let $|G|=2^{2k+3}$. Then
$g^2=t_1t_4^{1+r}(x_1^{-2}x_2)^{s-r}$ and so $g$ has order $2$ if and only
if $t_1t_4^{1+r}(x_1^{-2}x_2)^{s-r}=1,$ $0\le s,r\le 2^k-1.$ For arbitrary values of
$t_1$ and $t_4$ the last equality is valid for exactly $2^k$ different pairs $(s,r).$
Now notice that $(yx_1^{-1})^{y^2x_1^{-2}}=(yx_1^{-1})t_1t_4$ and
$(yx_1^{-1}x_2^{2^{k-1}})^{y^2x_1^{-2}x_2^{2^{k-1}}}=(yx_1^{-1}x_2^{2^{k-1}})zt_1t_4.$
Hence $t_1=t_4$ implies that $yx_1^{-1}$ has order $2,$
$C_G(yx_1^{-1})=\langle yx_1^{-1},\ y^2x_1^{-2},\
x_1^{-2}x_2\rangle $ and $|C_{yx_1^{-1}}|=2^{k+1}.$ Therefore all
elements of order $2$ of the set $M_3\setminus H$ lie in
$C_{yx_1^{-1}}.$ If $t_1\neq t_4,$ then $yx_1^{-1}x_2^{2^{k-1}}$
has order $2$  and similarly as in the previous
case, all elements of order $2$ of the set $M_3\setminus H$ lie in
the conjugacy class represented by $yx_1^{-1}x_2^{2^{k-1}}$ since
$C_G(yx_1^{-1}x_2^{2^{k-1}})=\langle yx_1^{-1}x_2^{2^{k-1}},\
y^2x_1^{-2}x_2^{2^{k-1}},\ x_1^{-2}x_2\rangle .$
\end{proof}

\begin{corollary}
Tables 15,16 and 17 contain full information about orders of elements of $X$ in
particular groups.
$$
\begin{array}{|c|c|c|c|c|c|c|c|c|}
\hline
g                      & G_{28} & G_{29} & G_{30} & G_{31} & G_{32} & G_{33} & G_{34} & G_{35} \\
\hline
y^2                    &    2   &    2   &    4   &    4   &    2   &    2   &    4   &    4   \\
y^2x_1^{-2}            &    2   &    2   &    2   &    2   &    4   &    4   &    4   &    4   \\
yx_1^{-1}              &    2   &    4   &    2   &    4   &    2   &    4   &    2   &    4   \\
yx_1^{-1}x_2^{2^{k-1}} &    2   &    4   &    2   &    4   &    2   &    4   &    2   &    4   \\
\hline
\end{array}\\
$$
\begin{center}
{\rm Table 15.}
\end{center}
$$
\begin{array}{|c|c|c|c|c|}
\hline
g                      & G_{36} & G_{37} & G_{38} & G_{39} \\
\hline
y^2                    &    2   &    2   &    4   &    4   \\
yx_1^{-1}x_2^{2^{k-2}} &    2   &    4   &    2   &    4   \\
\hline
\end{array}
\ \ \begin{array}{|c|c|c|c|c|}
\hline
g                     &  G_{40} & G_{41} & G_{42} & G_{43}\\
\hline
y^2                   &     2   &    4   &    2   & 4  \\
y^2x_1^{-2}           &     2   &    2   &    4   & 4  \\
yx_1^{-1}             &     2   &    2   &    4   & 4  \\
yx_1^{-1}x_2^{2^{k-1}} &     4   &    4   &    2   & 2  \\
\hline
\end{array}
$$
\hspace{3cm} {\rm Table 16. \hspace{4cm} Table 17.}
\end{corollary}

\vspace{8pt}
Now, similarly as in the case of $Fam8$ we can use it for building maximal
elementary abelian subgroups and counting the Quillen parameters  for our groups.

\begin{prop}
Let $G=G_m\in Fam7$. Then the Quillen parameter of $G$ is such as in the following table
$$
\begin{array}{|c||c|l|}\hline
{\rm Group} & {\rm Quillen}   & \text{Representatives of the conjugacy classes} \\
            & {\rm parameter} & \text{of maximal elementary abelian subgroups} \\ \hline
G_{28}      & (0,0,1,1)       & \langle y^2, \Omega_1(A) \rangle,\ \langle yx_1^{-1}, y^2x_1^{-2}, \Omega(A)\rangle  \\ \hline
G_{29}      & (0,0,2,0)       & \langle y^2, \Omega_1(A) \rangle,\ \langle y^2x_1^{-2}, \Omega_1(A)\rangle \\ \hline
G_{30}      & (0,0,0,1)       & \langle yx_1^{-1}, y^2x_1^{-2}, \Omega_1(A)\rangle \\ \hline
G_{31}      & (0,0,1,0)       & \langle y^2x_1^{-2},\Omega_1(A) \rangle \\ \hline
G_{32}      & (0,0,2,0)       & \langle y^2, \Omega_1(A) \rangle,\ \langle yx_1^{-1}, \Omega_1(A)\rangle  \\ \hline
G_{33}      & (0,0,1,0)       & \langle y^2, \Omega_1(A)\rangle \\ \hline
G_{34}      & (0,0,1,0)       & \langle yx_1^{-1},\Omega_1(A)\rangle \\ \hline
G_{35}      & (0,1,0,0)       & \Omega_1(A) \\ \hline
G_{36}      & (0,0,2,0)       & \langle y^2, \Omega_1(A) \rangle,\ \langle yx_1^{-1}x_2^{2^{k-2}}, \Omega_1(A)\rangle  \\ \hline
G_{37}      & (0,0,1,0)       & \langle y^2, \Omega_1(A) \rangle \\ \hline
G_{38}      & (0,0,1,0)       & \langle yx_1^{-1}x_2^{k-2}, \Omega_1(A)\rangle \\ \hline
G_{39}      & (0,1,0,0)       & \Omega_1(A) \\ \hline\hline
G_{40}      & (0,0,3,0)       & \langle y^2, \Omega_1(A) \rangle, \langle y^2x_1^{-2}, \Omega_1(A)\rangle, \langle yx_1^{-1}, y^2x_1^{-2},z \rangle  \\ \hline
G_{41}      & (0,0,2,0)       & \langle y^2x_1^{-2}, \Omega_1(A)\rangle, \langle yx_1^{-1}, y^2x_1^{-2},z \rangle \\ \hline
G_{42}      & (0,1,1,0)       & \langle y^2, \Omega_1(A)\rangle, \langle yx_1^{-1}x_2^{2^{k-1}}, z \rangle \\ \hline
G_{43}      & (0,2,0,0)       & \Omega_1(A), \langle  yx_1^{-1}x_2^{2^{k-1}}, z \rangle \\ \hline
\end{array}\label{fam7-Q}
$$
\begin{center}
{\rm Table 18.}
\end{center}
\end{prop}

For counting the numbers of conjugacy classes and the Roggenkamp parameters
we need much more detailed considerations than in the previous families.

\vspace{8pt}
\begin{lemma}\label{fam7-A-classes}
Let $G\in Fam7.$\\
(a) If $|G|=2^{2k+3},$ then $|Cl_G(A)|=2^{2k-3}+2^{k-1}+2^{k-2}+1;$\\
(b) If $|G|=2^{2k+2}$ and $A$ is abelian, then $|Cl_G(A)|=2^{2k-4}+2^{k-1}+1$;\\
(c) If $|G|=2^{2k+2}$ and $A$ is nonabelian, then
$$|Cl_G(A)|=2^{2k-5}+2^{2k-7}+
2^{k-2}+2^{k-3}+2^{k-4}+1.$$ \\
Moreover,
 $$R_G(A)=\left\{\begin{array}{lll}
 2|Cl_G(A)|   & {\rm if} & |G|=2^{2k+3}, \\
 2|Cl_G(A)|+1 & {\rm if} & |G|=2^{2k+2}.
 \end{array}\right.$$
\end{lemma}

\begin{proof}
(a) Since $A$ is abelian and $|G:A|=8,$ each conjugacy class of
$G$ contained in $A$ has at most $8$ elements. For a noncentral
element $g$ of $\Omega_1(A)=\langle (x_1^2)^{2^{k-1}},
x_2^{2^{k-1}}\rangle$ we have $C_G(g)=\langle y^2, x_1,
x_2\rangle=\langle y^2, x_1\rangle$ and then $|C_g|=2.$ For $g\in
A\setminus\Omega_1(A),$ $g^{y^2}=g^{-1}h\neq g,$ where $h\in
\{1,z\}$. Hence $y^2\notin C_G(g)$ and because of that
\mbox{$|G:C_G(g)|\geq 4,$} which says that either $C_G(g)=A$ or
$C_G(g)=H_i$ for certain $i\in\{1,2,3,4\}.$ Thus for $g\in
A\setminus \Omega_1(A)$, $|C_g|\le 4$ if and only if $g\in
T=C_A(x_1)\cup C_A(y^2x_1)\cup C_A(yx_1)\cup C_G(y^3x_1).$
Straightforward computations show that $C_A(x_1) = \langle
x_1^2,\Omega_1(A)\rangle,\  C_A(y^2x_1) = \langle
x_2,\Omega_1(A)\rangle, $ $C_A(yx_1) = \langle
x_1^2x_2^{-1}\rangle,\ C_A(y^3x_1) = \langle x_1^2x_2\rangle.$ It
is clear that $C_A(x_1)\cap C_A(y^2x_1)=\Omega_1(A)$ and
$C_A(yx_1)\cap C_A(y^3x_1)=\langle (x_1^2x_2)^{2^{k-1}}\rangle
\le\Omega_1(A).$ Hence $|C_A(x_1)\cup C_A(y^2x_1)|=2^{k+2}-4,$
$|C_A(yx_1)\cup C_A(y^3x_1)|=2^{k+1}-2$ and so
$|T|=2^{k+2}+2^{k+1}-8.$ The set $T$ is of course a normal subset
of G and it contains three conjugacy classes ha\-ving less than $4$
elements, namely $\{e\},\ \{(x_1^2x_2)^{2^{k-1}}\}$ and
$\{(x_1^2)^{2^{k-1}}, x_2^{2^{k-1}}\}$. Therefore
$|Cl_G(T)|=3+\frac{2^{k+2}+2^{k+1}-12}{4}=2^{k}+2^{k-1}$. All
other elements of $A$ lie in $8$-element conjugacy classes. So we
have
$\frac{2^{2k}-2^{k+2}-2^{k+1}+8}{8}=2^{2k-3}-2^{k-1}-2^{k-2}+1$
such classes. Consequently
$|Cl_G(A)|=2^{2k-3}+2^{k}-2^{k-2}+1=2^{2k-3}+2^{k-1}+2^{k-2}+1.$

(b) Similarly as in the proof of the first case for $a\in A$,
$|C_a|=4$ if and only if $C_G(a)$ is one of the five subgroups of
index $4$ containing $A.$ If $a\notin \Omega_1(A),$ then
$a^{y^2}=a^{-1}t,$ where $t\in\Omega_1(A)$, so $H$ is not a
centralizer of an element from $A\setminus \Omega_1(A).$ If
$C_G(a)=H_1,$ then $a\in\langle x_1^2\rangle =A_1.$ Since
$(x_1^2)^y=x_2t_3t_4$, $C_G(a)=H_2$ for $a\in \langle
x_2t_3t_4\rangle=A_2.$ Further, if $C_G(a)=H_3,$ then $a\in
\langle x_1^{-2}x_2,\Omega_1(A) \rangle =B_1$ and similarly
$C_G(a)=H_4,$ if $a\in \langle x_1x_2, \Omega_1(A) \rangle=B_2.$
Since $A_1A_2=A$ and $|A_1|=|A_2|=2^k$ we obtain $|A_1\cap
A_2|=2.$ Analogously, $|A:B_1B_2|=2$ and $|B_1|=|B_2|=2^k,$ so
$|B_1\cap B_2|=4.$ Hence $|A_1\cup A_2|=2^{k+1}-2,$ $|B_1\cup
B_2|=2^{k+1}-4$ and from $|(A_1\cup A_2)\cap (B_1\cup B_2)|=2$ we
obtain $|A_1\cup A_2\cup B_1\cup B_2|= 2^{k+2}-8.$ Since $A_1\cup
A_2\cup B_1\cup B_2$ contains all elements, whose conjugacy
classes have no more than $4$ elements, we obtain
$3+\frac{2^{k+2}-12}{4}=2^k$ $4$-element classes. Each of all
other classes contained in $A$ has $8$ elements. So we have
$\frac{2^{2k-1}-2^{k+2}+8}{8}=2^{2k-4}-2^{k-1}+1$ such classes.
Finally we obtain
$|Cl_G(A)|=2^{2k-4}+2^k-2^{k-1}+1=2^{2k-4}+2^{k-1}+1.$\par (c) It
is clear that $A$ is nonabelian only when $x_2^{x_1}=x_2^{-1}z_2,$
that is when $G\in\{G_{36},\ G_{37},\ G_{38},\ G_{39}\}.$ In this
case the action of $G$ on $A^2$ is similar to the action of $G$ on
$A$ in the previous case, that is when $A$ is abelian and
$|A|=2^{2k-3}.$ So $|Cl_G(A^2)|=2^{2k-6}+ 2^{k-2}+1.$ Conjugacy
classes which are contained in $A\setminus A^2$ have either $8$ or
$16$ elements. Since for every $g\in A\setminus A^2,$
$C_g=C_{gz_1}$ and $\zeta(G)= \langle z_1\rangle,$ images of these
classes in $\widetilde{G}=G/\zeta(G)$ have twice less elements and
images of different classes are different. By (a)
$|Cl_{\widetilde{G}}(\widetilde{A}\setminus
\widetilde{A}^2)|=(2^{2k-5}+2^{k-2}+2^{k-3}+1)-(2^{2k-7}+2^{k-3}+2^{k-4}+1)=
2^{2k-5}-2^{2k-7}+2^{k-3}+2^{k-4}.$ Therefore $|Cl_{G}(A\setminus
A^2)|=2^{2k-5}-2^{2k-7}+2^{k-3}+2^{k-4}+2^{2k-6}+2^{k-2}+1=
2^{2k-5}+2^{2k-7}+2^{k-2}+2^{k-3}+2^{k-4}+1.$\par For the proof of
the last assertion it is enough to notice that $d(M_2)=2$,
$d(M_3)=3$ and if $g\in \Omega_1(A)\setminus \zeta(G),$ then
$$C_G(g)=\begin{cases}
M_2 & \text{if $|G|=2^{2k+3}$,} \\
M_3 & \text{if $|G|=2^{2k+2}$.}
\end{cases}$$
All other elements of $A$ have centralizers equal to $H_i$ for
some $i,$ $1\leq i\leq 4,$ or to $G$. All these groups are
$2$-generated.
\end{proof}

\begin{lemma} \label{fam7-M2}
Let $G\in Fam7$ and let $x\in M_2\setminus H.$\par
(a) If
$|G|=2^{2k+3},$ then $C_G(x)= \langle x,\Omega_1(A)\rangle $ and
$|C_x|=2^{k+1}.$ Moreover, $$|Cl_G(M_2\setminus H)|=2^k \ \ {\rm
and}\ \ R_G(M_2\setminus H)= 2^{k+1}.$$\par
(b) If $|G|=2^{2k+2},$
then $C_G(x)= \langle x\rangle $ and $|C_x|=2^{k+1}.$ Moreover,
$$|Cl_G(M_2\setminus H)|=2^{k-1}=R_G(M_2\setminus H).$$
\end{lemma}

\begin{proof}
Since $M_2\setminus H=x_1A\cup y^2x_1A$ and $(x_1A)^y=y^2x_1A,$ it
suffices to count $|C_x|$ for $x\in x_1A.$ Put
$x=x_1x_1^{2r}x_2^s\in  x_1A.$ If $g\in
G$ centralizes $x$ then $\overline{g}=gA$ centralizes
$\overline{x}=xA$ in $\overline{G}=G/A$ that is $\overline{g}\in
\langle \overline{y}^2, \overline{x}_1\rangle.$ Now, it follows
from properties of $\overline{G}$ that $g=y^2a$ or $g=y^2x_1a$ or
$g\in \langle x_1,A\rangle $ for suitable $a=x_1^{2u}x_2^{v}\in
A.$ None of the elements of the form either $g=y^2a$ or $g=y^2x_1a$
centralizes $x.$ In fact,
$$x^{y^2a}=(x_1x_1^{2r}x_2^s)^{y^2x_1^{2u}x_2^v}=
x_1x_1^{-2r+2}x_2^{-s+1+2v}a_1,$$ where $a_1\in \Omega_1(A).$
Similarly, $$x^{y^2x_1a}=(x_1x_1^{2r}x_2^s)^{y^2x_1x_1^{2u}x_2^v}=
x_1x_1^{-2r+2}x_2^{s-1+2v}a_1.$$ So in both cases we obtain $x\neq
x^g.$ Therefore we assume that $g\in \langle x_1,A\rangle .$ In
this case we need only to describe $C_A(x).$ For $x=x_1$ it was
done in the proof of Lemma \ref{fam7-A-classes}. For the general
case the calculations and conclusions are similar. So for a) we
obtain $C_G(x)=\langle x,\Omega_1(A)\rangle$ and $d(C_G(x))=2.$
Hence $|C_x|=\frac{2^{2k+3}}{2^{k+2}}=2^{k+1}.$ This means that
$|Cl_G(M_2\setminus H)|=\frac{|M_2\setminus H|}{2^{k+1}}=2^k$ and
$R_G(M_2\setminus H)=2|Cl_G(M_2\setminus H)|=2^{k+1}.$ \par For b)
the proof of Lemma \ref{fam7-A-classes} yields $C_G(x)=\langle
x\rangle$ (and of course $d(C_G(x))=1$). Thus $|Cl_G(M_2 \setminus
H)|=\frac{|M_2\setminus H|}{2^{k+1}}= 2^{k-1} = R_G(M_2 \setminus
H).$
\end{proof}

For counting the centralizers of elements of $M_3\setminus H$
notice that $M_3\setminus H=yx_1^{-1}A\cup y^3x_1A$,
$(yx_1^{-1}A)^{x_1}=y^3x_1A$ and $(y^3x_1A)^{x_1}=yx_1^{-1}A$. So
it suffices to study only these conjugacy classes which are
represented by elements of the form $g=yx_1^{-1}a$, where
$a=x_1^{2r}x_2^s\in A$.

\begin{lemma}\label{fam7-M3-1}
Let $G\in Fam7$ and let $A$ be abelian. If $g\in yx_1^{-1}A$, then
$$C_A(yx_1^{-1}a)=
\left\{\begin{array}{lll}
\langle x_1^2x_2^{-1}\rangle      & {\rm if} & |G|=2^{2k+3}, \\
\langle x_1^2x_2^{-1}, z_1\rangle & {\rm if} & |G|=2^{2k+2}.
\end{array}\right. $$
\end{lemma}

\begin{proof}
For every $h=x_1^{2u}x_2^w\in A$ we obtain from the relations
of the group and the relations (\ref{eq1}) that

$$\begin{tabular}{rl} \vspace{12pt} $y^{h}=$ &
$(x_2^{-w}x_1^{-2u})y(x_1^{2u}x_2^w)
=yx_1^{2w}x_2^{-u}t_4^{-u}x_1^{2u}x_2^w=
yx_1^{2(u+w)}x_2^{-u+w}t_4^{-u}$\\
\vspace{6pt} $(x_1^{-1})^h=$
         & $x_2^{-w}x_1^{-1}x_2^w= x_1^{-1}x_2^{2w}t_4^{-w}$ \\
\end{tabular}$$

\noindent Therefore

\begin{equation}\label{eq-fam7-M3-1}
\vspace{6pt}\begin{tabular}{rl}
\vspace{6pt}
$(yx_1^{-1})^{h}=$ & $(yx_1^{2(u+w)}x_2^{-u+w}t_4^{-u})(x_1^{-1}x_2^{2w}t_4^{-w})=$\\
\vspace{6pt}
         & $yx_1^{-1}x_1^{2(u+w)}x_2^{u-w}t_4^{-2u+w}x_2^{2w}t_4^{-w})=$ \\
\vspace{6pt}
         & $yx_1^{-1}(x_1^2x_2)^{u+w}$
\end{tabular}
\end{equation}

\noindent Hence for every $a=x_1^{2r}x_2^s\in A$,
$h=x_1^{2u}x_2^w$ centralizes $yx_1^{-1}a$ if and only if
$(x_1^2x_2)^{u+w}=1$. This means for the case $|G|=2^{2k+3}$ that
$u+w=0({\rm mod}\, 2^k)$ Consequently $h=(x_1^{-2}x_2)^w$ and then
$C_A(yx_1^{-1}a)=\langle (x_1^{-2}x_2)\rangle$. If $|G|=2^{2k+2}$,
then we get $u+w=0({\rm mod}\, 2^{k-1})$ or equivalently $u+w=0\ {\rm
or}\ 2^{k-1}({\rm mod}\, 2^k)$. Hence $h=(x_1^{-2}x_2)^w$ or
$h=(x_1^{-2}x_2)^wz_1$. Finally $C_A(yx_1^{-1}a)=\langle
x_1^{-2}x_2,z_1\rangle$.
\end{proof}

The centralizer of the image of $g=yx_1^{-1}a$ in the factor group
$\overline{G}=G/A$ is equal to $\langle \overline{g},
\overline{y}^2\rangle.$ So we need also to check for which $g\in
yx_1^{-1}A$ there exist elements in $y^2A$ centralizing $g.$

\begin{lemma}\label{fam7-M3-2}
Let $G=G_n\in Fam7$ and let $A$ be abelian. \par (a) If
$n\in\{28,29,30,31\}$ and $g\in \{yx_1^{-1},
yx_1^{-1}x_2^{2^{k-1}}\}$ or $n\in\{32,33,34,35\}$ and
$g\in\{yx_1^{-1}x_2^{2^{k-2}}, yx_1^{-1}x_2^{-2^{k-2}}\}$, then
$C_G(g)=\langle
g,y^2x_2^{-1},x_1^{-2}x_2,x_2^{2^{k-1}}\rangle$.\par (b) If
$n\in\{40,41,42,43\}$ and $g\in \{yx_1^{-1},
yx_1^{-1}x_2^{2^{k-1}}\}$, then $C_G(g)=\langle g, y^2x_2^{-1},$\par\noindent
$x_{1}^{-2}x_2\rangle$.\par (c) If $g\in yx_1^{-1}A$ is not
conjugated to elements distinguished in (a) and (b) then
$C_G(g)=\langle g,x_1^{-2}x_2, \zeta(G) \rangle$.\par (d)
$Cl_G(M_3\setminus A)=2^{k-1}+1.$
\end{lemma}

\begin{proof}
It follows from Lemma \ref{fam7-M3-1} that in all the cases
$|C_A(g)|=2^{k}$ for all $g\in yx_1^{-1}A.$ It is also easily seen
that $g^2\in C_A(g)$ and for every $h\in A$, $(y^2h)^2\in C_A(g).$
Thus if $C_G(g)=\langle g,y^2h,C_A(g)\rangle,$ then
$|C_G(g)|=2^{k+2}$ and so
$$|C_g|=\begin{cases}
2^k     & \text{for $|G|=2^{2k+2}$,}\\
2^{k+1} & \text{for $|G|=2^{2k+3}$.}
\end{cases}$$
All other conjugacy classes have $2^{k+1}$ and $2^{k+2}$ elements
in the cases $|G|=2^{2k+2}$ and $|G|=2^{2k+3}$ respectively.

Now let $h=x_1^{2u}x_2^w$ and $a=x_1^{2r}x_2^{s}$ be arbitrary
elements of $A$. As it was seen in the proof of Lemma
\ref{fam7-M3-1}
$$y^{y^2h}=y^h=yx_1^{2(u+w)}x_2^{-u+w}t_4^{-u}.$$ Moreover

\vspace{6pt}
\begin{tabular}{rl}
\vspace{6pt}$(x_1^{-1})^{y^2h}=$
       & $(x_1^{-1})^{y^2x_1^{2u}x_2^{w}}=
       (t_1x_2^{-1}x_1)^{x_1^{2u}x_2^{w}}=
       (t_1x_2^{-1}x_1)^{x_2^{w}}=$\\
\vspace{6pt}
       & $x_2^{-w}(t_1x_2^{-1}x_1)x_2^{w}=
       x_1t_4^{-w}x_2^{w}(t_1t_4^{-1}x_2)x_2^{w}=$\\
\vspace{6pt}
       & $x_1x_2^{2w+1}(t_1t_4^{-1-w})$\\
\end{tabular}

\noindent and
$$a^{y^2h}=(x_1^{2r}x_2^{s})^{y^2}=x_1^{-2r}x_2^{-s}t_4^{r+s}=
a^{-1}t_4^{r+s}.$$
Hence
\begin{equation}\label{eq-fam7-M3-2}
\begin{split}
(yx_1^{-1}a)^{y^2x_1^{2u}x_2^w}= &
(yx_1^{2(u+w)}x_2^{-u+w}t_4^{-u})(x_1x_2^{2w+1}(t_1t_4^{-1-w}))
       (a^{-1}t_4^{r+s})=\\
& yx_1^{-1 }(x_1^2x_2)^{u+w+1}t_4^{-1}(a^{-1}t_4^{r+s})= \\
& yx_1^{-1}a(x_1^2x_2)^{u+w+1}a^{-2}(t_1t_4^{r+s-1}).
\end{split}
\end{equation}
Since $t_1,t_4\in \zeta(G)$, we obtain that if $y^2h\in
C_G(yx_1^{-1}a)$, then $a^2\in \langle x_1^2x_2, \zeta(G) \rangle$. If
$|G|=2^{2k+3},$ then the subgroup $\langle x_1^2x_2, \zeta(G) \rangle$
is cyclic of order $2^k$. Therefore for every $a\in A$ satisfying
$a^2\in \langle x_1^2x_2, \zeta(G) \rangle$ we can find $u$ and $w$
such that $y^2(x_1^{2u}x_2^w)\in C_G(yx_1^{-1}a)$. There exists
exactly $2^{k+1}$ possible values for $a$, because $a$ can be an
arbitrary element of $\langle x_1^2x_2, \Omega_1(A)\rangle.$ It is
an easy task to show that for $g=yx_1^{-1},$
$$\begin{tabular}{rl}
$C_g=$ & $ \{yx_1^{-1}(x_1^{2}x_2)^r: \ r=0,1,\ldots
,2^{k}-1\}\cup$ \\ & $\{y^3x_1^{-1}(x_1^{-2}x_2)^r: \ r=0,1,\ldots
,2^{k}-1\} $\end{tabular}$$
and for $g=yx_1^{-1}x_2^{2^{k-1}},$

$$\begin{tabular}{rl}
$C_g=$ & $\{yx_1^{-1}(x_1^2x_2)^rx_2^{2^{k-1}}: \ r=0,1,\ldots
,2^{k}-1\}\cup$ \\ & $\{y^3x_1^{-1}(x_1^{-2}x_2)^rx_2^{2^{k-1}}: \
r=0,1,\ldots ,2^{k}-1\}.$ \end{tabular}$$

Both classes contain obviously all
elements of $yx_1^{-1}A$ having centralizer of order $2^{k+2}$.
All other elements $g\in M_3\setminus H$ have centralizers of the
form $\langle g,x_1^{-2}x_2\rangle,$ because no element from the
coset $y^2A$ centralizes $g$.

Now assume that $|G|=2^{2k+2}$. Then $\langle
x_1^{2}x_2, \zeta(G) \rangle=\langle x_1^{2}x_2,z_1\rangle$ is not
cyclic and if for $a=x_1^{2r}x_2^s$ we have $a^2\in \langle
x_1^{2}x_2,z_1\rangle$, then $a\in \langle
x_1^{2}x_2,\Omega_2(A)\rangle$. Since we assumed in the beginning
of this section that $|G|\geq 2^7$ we have $k>2$ and then in
(\ref{eq-fam7-M3-2}) we obtain $ r+s \equiv 0({\rm mod}\ 2)$. If
$t_1t_4=1$ then $y^2h\in C_G(g)$, $g=yx_1^{-1}a$, if and only if
$a$ belongs to $\langle x_1^{2}x_2,\Omega_1(A)\rangle$ which is of
order $2^k.$ So there exist in $yx_1^{-1}A$ exactly $2^{k}$
elements with the centralizers of order $2^{k+2}.$ It is clear
that they are divided into two conjugacy classes as for such $g$
we have $|C_g\cap yx_1^{-1}A|=2^{k-1}.$ It is also clear that the
elements $yx_1^{-1}$ and $yx_1^{-1}z_1$ are not conjugated and if
$g$ is one of them then $C_G(g)=\langle g,
y^2x_2^{-1},\Omega_1(A)\rangle.$ All other elements have the
centralizers of the form $\langle g,C_A(g)\rangle$ that is of
order $2^{k+1}$.

If $t_1t_4=z_1$, then $y^2h\in C_G(g)$, $g=yx_1^{-1}a$, if and
only if $a$ belongs to $\langle x_1^{2}x_2,\Omega_2(A)\rangle$ and
$a^2=z_1$. So there exist again $2^{k}$ possible values for $a$.
Let $a=x_2^{\pm 2^{k-1}}$, $u=0 $ and $w=-1$. Then from
(\ref{eq-fam7-M3-2}) we obtain

$$
\begin{tabular}{rl}
\vspace{6pt}$(yx_1^{-1}x^{\pm 2^{k-2}})^{y^2x_2^{-1}}=$
       & $yx_1^{-1}x_2^{\pm 2^{k-2}}x_2^{-2^{k-1}}z_1=
       yx_1^{-1}x_2^{\pm 2^{k-2}}$
\end{tabular}
$$

Since the elements $yx_1^{-1}x^{2^{k-2}}$ and
$yx_1^{-1}x^{-2^{k-2}}$ are not conjugated the conjugacy classes
of them contain all the elements with the centralizers of order
$2^{k+2}.$ All other elements have the centralizers of the form
$\langle g,C_A(g)\rangle$ that is of order $2^{k+1}$.

For the proof of (d) let us notice that if $|G|=2^{2k+3}$ then we
have $2$ conjugacy classes having $2^{k+1}$ elements. All other
classes have twice more elements. Since $|M_3\setminus
H|=2^{2k+1}$ we have
$2+\frac{2^{2k+1}-2^{k+2}}{2^{k+2}}=2^{k-1}+1$ conjugacy classes.
The calculation for the case $|G|=2^{2k+2}$ is similar.

\end{proof}

\begin{lemma} \label{fam7-M3-R}
Let $G=G_n\in Fam7$, $|G|=2^{2k+2+\epsilon}$, where
$\epsilon\in\{0,1\}$. If $A$ is abelian, then
$$R_G(M_3\setminus H)=\left\{\begin{array}{ll}
2^{k}+2^{k-2}+5 & {\rm if}\ n\in\{28,30\} \\
2^{k}+2^{k-2}+3 & {\rm if}\ n\in\{32,34\} \\
2^{k}+4           & {\rm if}\ n\in\{29,31,33,35\} \\
2^{k-1}+2^{k-2}+4 & {\rm if}\ n\in\{40,41,42,43\} \\
\end{array}\right.$$

\end{lemma}

\begin{proof}
First assume that $|G|=2^{2k+3}$ i.e. $n\in\{40,41,42,43\}.$ If
$g=yx_1^{-1}$ or $g=yx_1^{-1}x_2^{2^{k-1}}$ then $C_G(g)= \langle
g,y^2x_1^{-2},x_1^{-2}x_2\rangle $ by Lemma \ref{fam7-M3-2}(b).
Thus $d(C_G(g))=3$. Hence it remains to consider conjugacy classes
contained in $M_3\setminus H$ containing neither $yx_1^{-1}$ nor
$yx_1^{-1}x_2^{2^{k-1}}.$ Let $g$ be a representative of such a
class. It follows from Lemma \ref{fam7-M3-2}(c) that
$C_G(g)=\langle g,x_1^2x_2\rangle$. Thus this centralizer is
either $2$-generated or cyclic. Put
$g=yx_1^{-1}(x_1^{2r}x_2^{s}).$ So
$g^2=(x_1^2x_2^{-1})^{r-s}t_1t_4$ and then $C_G(g)$ is cyclic if
and only if $r-s$ is invertible in $\mathbb{Z}_{2^k}.$ There are
exactly $2^{2k-1}$ such elements. Since they split into conjugacy
classes, each having $2^{k+1}$ elements in $yx_1^{-1}A$, we have
$2^{k-2}$ such classes that their representatives have cyclic
centralizers. Representatives of other classes contained in
$M_3\setminus H$ have $2$-generated centralizers. Therefore by
Lemma \ref{fam7-M3-2}(d) $$R_G(M_3\setminus H)=6+ 2^{k-2}+ 2\cdot
(2^{k-2}-1)= 2^{k-1}+2^{k-2}+4.$$

Now assume that $|G|=2^{2k+2}$. If $g$ is one of the elements
distinguished in Lemma \ref{fam7-M3-2}(a) then $C_G(g)=\langle g,
y^2x_2^{-1},x_1^2x_2,z_1\rangle $ and then $d(C_G(g))=4$ if and
only if $o(g)=o(y^2x_2^{-1})=2$. It follows from Table $12$ that
this happens for the groups $G_{28}$ and $G_{30}$ only. In all
other cases we have $d(C_G(g))=3$. So it remains to consider
conjugacy classes contained in $M_3\setminus H$ not containing
elements distinguished in Lemma \ref{fam7-M3-2}(a). It follows
from Lemma \ref{fam7-M3-2}(c) that $C_G(g)=\langle
g,x_1^{-2}x_2,z_1\rangle$. Thus this centralizer is either
$2$-generated or $3$-generated. Put
$g=yx_1^{-1}(x_1^{2r}x_2^{s}).$ So
$g^2=t_1t_2t_4(x_1^{-2}x_2)^{s-r},$ and then $C_G(g)$ is
$2$-generated if and only if $t_1t_2t_4=z_1$ or $t_1t_2t_4=1$ and
$s-r$ is invertible in $Z_{2^{k-1}}.$ By Table 5 $t_1t_2t_4=z_1$
in $G_m$ if $m\in \{29,31,33,35\}$ only. So in these groups
$d(C_G(g))=2.$ If $m\in \{28,30,32,34\}$, then $t_1t_2t_4=1$ and
we have exactly $2^{2k-2}$ such elements $g$ that $r-s$ is
invertible. Since they split into conjugacy classes having $2^{k}$
elements from $yx_1^{-1}A$ we have $2^{k-2}$ such classes that
their representatives have $2$-generated centralizers.
Representatives of other classes contained in $M_3\setminus H$
have $3$-generated centralizers. Therefore the Roggenkamp number
of $M_3\setminus H$ is equal to

$$R_G(M_3\setminus H)=\left\{
\begin{array}{l}
2\cdot 2^{k-2}+ 3\cdot (2^{k-2}-1)+8\\
2\cdot 2^{k-2}+ 3\cdot (2^{k-2}-1)+6\\
2\cdot (2^{k-1}-1)+6
\end{array}\right.=$$

$$= \left\{
\begin{array}{ll}
2^k +2^{k-2}+5 & {\rm if}\ \ m\in\{28,30\},\\
2^k +2^{k-2}+3 & {\rm if}\ \ m\in\{32,34\},\\
2^k+4          & {\rm if}\ \ m\in\{29,31,33,35\}.
\end{array}
\right.$$

\end{proof}

\begin{prop}\label{fam7-odd-Cl}
If $G$ is a group of $Fam7$ with $A$ abelian then
$$|Cl(G)|=\left\{\begin{array}{ll}
2^{2k-4}+3\cdot 2^{k-1}+6 & {\rm if }\ |G|=2^{2k+2},\\
2^{2k-3}+9\cdot 2^{k-2}+6 & {\rm if }\ |G|=2^{2k+3}.

\end{array}\right. $$

\end{prop}

\begin{proof}
For the proof it is enough to count the conjugacy classes
described in lemmas \ref{fam7-M1}--\ref{fam7-M3-2}. We gather this
information in the following table.
$$\begin{array}{|c|c|c|}
\hline {\rm Subset}   & |G|=2^{2k+2}         & |G|=2^{2k+3}\\
\hline A       & 2^{2k-4}+2^{k-1}+1   & 2^{2k-3}+2^{k-1}+2^{k-2}+1 \\
\hline H\setminus A   & 2             & 2     \\
\hline M_1\setminus H & 2             & 2     \\
\hline M_2\setminus H & 2^{k-1}       & 2^{k} \\
\hline M_3\setminus H & 2^{k-1}+1     & 2^{k-1}+1\\
\hline \sum           & 2^{2k-4}+2^{k}+2^{k-1}+6 & 2^{2k-3}+2^{k+1}+2^{k-2}+6 \\ \hline
\end{array}$$
\begin{center}
{\rm Table 19.}
\end{center}

\end{proof}

\begin{prop}\label{fam7-odd-Rg}
Let $G=G_n$ be a group of $Fam7$ and let $A$ be abelian. Then
$$R(G)=\left\{\begin{array}{lll}
2^{2k-3}+11\cdot 2^{k-2}+r_n & {\rm if} & n\in\{28,30,32,34\};\\
2^{2k-3}+5\cdot 2^{k-1}+r_n & {\rm if} & n\in\{29,31,33,35\};\\
2^{2k-2}+17\cdot 2^{k-2}+r_n & {\rm if} & n\in\{40,41,42,43\}.\\
\end{array}\right.$$
where $r_n$ for all the groups is given in the following table:

$$\begin{array}{|c|c|c|c|c|c|c|c|c|c|c|c|}
\hline
G_{28} & G_{29} & G_{30} & G_{31} & G_{32} & G_{33} & G_{34} & G_{35} & G_{40} & G_{41} & G_{42} & G_{43} \\

\hline
18    &   16    &  16    &   14   &   14   &   15   &   12   &   13   &   15   &   13   &   15   &   13   \\
\hline
\end{array}$$
\begin{center}
{\rm Table 20.}
\end{center}

\end{prop}

\begin{proof}
For centralizers of representatives of conjugacy classes we need
to count numbers of generators in minimal generating sets. First
we list centralizers of representatives of conjugacy classes
contained in $M_2\setminus A=(M_2\setminus H)\cup(H\setminus A).$
The element $g$ appearing in the centralizer of $y^2x_2^{-1}$ is
one of the elements defined in Lemma \ref{fam7-M3-2}

$$
\begin{array}{|c|c|c|c|c|c|c|c|}\hline
 g                       & C_G(g)                                                                              & G_{28} & G_{29} & G_{30} & G_{31} & G_{32} & G_{33} \\ \hline\hline
 y^2                     & \langle y, \Omega_1(A) \rangle                                                      &    2   &    2   &   2    &    2   &    2   &  2 \\ \hline
 y^2x_2^{-1}             & \langle y^2x_2^{-1}, g,\Omega_1(A)\rangle                                   &    4   &    3   &   4    &    3   &    3   &  2 \\ \hline
 y                       & \langle y, z_1 \rangle                                                              &    2   &    2   &   1    &    1   &    2   &  2 \\ \hline
 y^3                     & \langle y, z_1 \rangle                                                              &    2   &    2   &   1    &    1   &    2   &  2 \\ \hline
                         & \sum                                                                                &   10   &    9   &   8    &    7   &    9   &  8 \\ \hline
\end{array}\label{fam7-even-ab}$$
\begin{center}
{Table 21.}
\end{center}

$$
\begin{array}{|c|c|c||c|c|}\hline
 g                       & C_G(g)                                                                              & G_{33},G_{34} & G_{40},G_{42} & G_{41}, G_{43} \\ \hline\hline
 y^2                     & \langle y, \Omega_1(A) \rangle                                                      &    2          &    2          &    2           \\ \hline
 y^2x_2^{-1}             & \langle y^2x_2^{-1}, g,\Omega_1(A)\rangle                                           &    2          &    3          &    3           \\ \hline
 y                       & \langle y, z_1 \rangle                                                              &    1          &    2          &    1           \\ \hline
 y^3                     & \langle y, z_1 \rangle                                                              &    1          &    2          &    1           \\ \hline
                         & \sum                                                                                &    6          &    9          &    7           \\ \hline
\end{array}\label{fam7-even-ab1}$$
\begin{center}
{\rm Table 22.}
\end{center}

Now the sum of numbers from the appropriate columns of the above
tables and formulas counted for $R_G(A)$, $R_G(M_2\setminus H)$
and $R_G(M_3\setminus H)$ in Lemmas \ref{fam7-A-classes},
\ref{fam7-M2} and \ref{fam7-M3-R} gives the formulas of the proposition.
\end{proof}

\begin{lemma} \label{fam7-M3-3} Let $|G|=2^{2k+2}$ and let
$A$ be nonabelian.\par (a) If either $g=yx_1^{-1}$ or
$g=yx_1^{-1}x_2^{2^{k-1}}$ then $C_G(g)=\langle g, x_1^2x_2^{-1},
y^2x_1^{-2}, \Omega_1(A) \rangle$ and $|C_g|=2^{k}.$\par (b) If
$g\in M_3\setminus H$ is not conjugated neither to $yx_1^{-1}$ nor
to $yx_1^{-1}x_2^{2^{k-1}}$ then $C_G(g)=\langle g,
x_1^2x_2^{-1},z_1 \rangle$ and $|C_g|=2^{k+1}.$\par (c)
$|Cl_G(M_3\setminus H)|=2^{k-2}+2^{k-3}+1.$\par
\end{lemma}

\begin{proof}
Because of similar reasons as in the proof of Lemma
\ref{fam7-M3-2} we study classes represented by elements of the
form $g=yx_1^{-1}x_1^{2r}x_2^s$ only. Since $C_G(g)/\zeta(G) \le
C_{G/\zeta(G)}(g\zeta(G))$ by Lemma \ref{fam7-M3-2}(b) there are no
elements centralizing $g$ in $(M_1\setminus H)\cup (M_2\setminus
H).$ Let $h=x_1^{2u}x_2^w\in A$ be an arbitrary element. Then by
straightforward computations
\begin{equation}\label{m3_1}
g^h=g(x_1^{2(u+w)}x_2^{u+w}z_1^{ur+ws+u}t_3^{u+w}).
\end{equation}
Hence $h\in C_G(g)$ if and only if
\begin{equation}\label{m3_2}
u+w\equiv 0({\rm mod}\, 2^k) \ \ {\rm or}\ \ u+w\equiv 2^{k-1}
({\rm mod}\, 2^k)
\end{equation}
and
\begin{equation}\label{m3_3}
ur+ws+u\equiv 0({\rm mod}\, 2).
\end{equation}
If r and s are of different parity then each solution of
(\ref{m3_2}) is also a solution of (\ref{m3_3}). In this case
there are no elements centralizing $g$ in $H\setminus A.$
Otherwise the image of $g$ in $\widetilde{G}=G/\zeta(G)$ would have
such a centralizer which is not possible by Lemma
\ref{fam7-M3-2}(b). The number of elements of $yx_1^{-1}A$ of the
form $g=yx_1^{-1}x_1^{2r}x_2^s$ such that $r$ and $s$ are of
different parity is equal to $2^{2k-2}.$ The centralizer of such
$g$ has order $2^{k+1}$, so $|C_g|=2^{k+1}.$ Half of elements of
$C_g$ lie in $yx_1^{-1}A$. Thus we have just
$\frac{2^{2k-2}}{2^k}=2^{k-2}$ of such classes. Note that in this
case $C_g\neq C_{gz_1}$ ($g$ and $gz_1$ are not conjugated), the
images of $C_g$ and $C_{gz_1}$ in $\widetilde{G}$ are equal and
have the same size as $C_g.$ Hence there are $2^{k-3}$ conjugacy
classes in $\widetilde{G}$ which are images of classes just
considered. In $\widetilde{M_3}\setminus\widetilde{A}$ there are
$2^{k-2}+1$ classes. So we have to fix images of which classes of
$M_3\setminus A$ are the remaining classes.

If $r$ and $s$ are of the same parity, then odd solutions of
(\ref{m3_2}) does not satisfy (\ref{m3_3}). Moreover in this case
$g$ is conjugated to $gz_1$ ($g^{x_1^{2}x_2}=gz_1$). So images in
$\widetilde{G}$ of different such classes are different but they
have twice less elements. Hence we have just
$2^{k-2}+1-2^{k-3}=2^{k-3}+1$ of such classes. Therefore
$Cl_G(M_3\setminus A)=2^{k-2}+2^{k-3}+1.$
\end{proof}

\begin{lemma} \label{fam7-M3-4} Let $|G|=2^{2k+2}$ and let
$A$ be nonabelian. If $k>3$, then $$R_G(M_3\setminus
H)=\left\{\begin{array}{lll}
2^{k-1}+2^{k-2}+2^{k-4}+4 & {\rm if} & G\in\{G_{36},G_{38}\}, \\
2^{k-1}+2^{k-2}+4         & {\rm if} & G\in\{G_{37},G_{39}\}.
\end{array}\right.$$ If $k=3,$ then for $G\in\{G_{36},G_{38}\}$,
$R_G(M_3\setminus H)=10$ and for $G\in\{G_{37},G_{39}\}$,
$R_G(M_3\setminus H)=8.$
\end{lemma}
\begin{proof}
We consider only the case $k>3$. We give only sketch arguments
because detail calculations are similar as in previous proofs.
First let $g$ be one of the elements $yx_1^{-1}$ and
$yx_1^{-1}x_2^{2^{k-2}}.$ The images of the classes represented by
these elements in $\widetilde{G}=G/\zeta(G)$ are different by Lemma
\ref{fam7-M3-2}(b) and have centralizers of the form $\langle
\widetilde{g},C_{\widetilde{A}}(\widetilde{g}),
\widetilde{y}^2\widetilde{h}\rangle,$ where $h\in A.$ So the
centralizer of $g$ is also of this form. If $n\in\{37,39\}$, then
both elements have $3$-generated centralizers. If $n\in\{36,38\}$,
then one element have $3$-generated centralizer and the other one
$4$-generated centralizer.

Now let $g=yx^{-1}x_1^{2r}x_2^s\in yx^{-1}A$ be an arbitrary
element with $r,s$ of different parity. Accordingly to the proof
of Lemma \ref{fam7-even-nab-Rg} there exist $2^{k-2}$ classes
represented by elements of this form. By (\ref{m3_1}),
$C_G(g)=\langle g,x_1^{-2u}x_2^u,z_1\rangle$, so by (\ref{nab-g2})
it is $2$-generated, since $x_1^{-2u}x_2^u\in\langle
g,z_1\rangle.$

We have still to consider the remaining $2^{k-3}-1$ classes. If
$n\in\{37,39\}$ then all representatives of these conjugacy
classes are $2$-generated. So assume that $n\in\{36,38\}.$ Let
$g=yx^{-1}x_1^{2r}x_2^s\in yx^{-1}A$ be an arbitrary element with
$r,s$ of the same parity. If $s-r$ is divisible by $4$ then
$C_G(g)=\langle g,x_1^{-4}x_2^2,z_1\rangle$ is $3$-generated.
There exist $2^{k-4}-1$ such classes. If $s-r$ is not divisible by
$4$ then $C_G(g)$ is $2$-generated and there exist $2^{k-4}$ such
classes.

 Summarizing, if $n\in\{37,39\}$ then $R_G(M_3\setminus
H)=6+2\cdot2^{k-2}+2\cdot (2^{k-3}-1)=2^{k-1}+2^{k-2}+4$.
If $n\in\{36,38\}$ then $R_G(M_3\setminus
H)=7+2\cdot2^{k-2}+3\cdot (2^{k-4}-1)+2\cdot
2^{k-4}=2^{k-1}+2^{k-2}+2^{k-4}+4$.

\end{proof}

\begin{prop}\label{fam7-even-nab-Cl}
Let $G$ be a group of $Fam7$, $|G|=2^{2k+2}.$ If $A$ is
nonabelian then $|Cl(G)|=5\cdot 2^{2k-7}+21\cdot 2^{k-4}+6.$
\end{prop}

\begin{proof}
It follows from Lemmas \ref{fam7-M1}, \ref{fam7-A-classes},
\ref{fam7-M2} and \ref{fam7-M3-3} that
$$|Cl(G)|=|Cl_G(A)|+|Cl_G(H\setminus A)|+\sum_{i=1}^3|Cl_G(M_i\setminus
H)|=$$ $$(2^{2k-5}+2^{2k-7}+2^{k-1}+2^{k-2}+2^{k-3}+2^{k-4}+1)+2+
(2+2^{k-1}+2^{k-2}+2^{k-3}+1)$$ which is equal to the formula
given in the proposition.
\end{proof}

\begin{prop}\label{fam7-even-nab-Rg}
Let $G=G_m$ be a group of $Fam7$, $|G|=2^{2k+2}.$ If $A$ is
nonabelian and $k>3$, then $$R(G_m)=\left\{\begin{array}{lll}
5\cdot 2^{2k-6}+35\cdot 2^{k-4}+r_m & {\rm if} & m\in\{36,38\}, \\
5\cdot 2^{2k-6}+17\cdot 2^{k-3}+r_m & {\rm if} & m\in\{37,39\},
\end{array}\right.$$
where $r_{36}=16,\ r_{37}=15,\ r_{38}=14$ and $r_{39}=13.$ For
$k=3$ we have $R(G_{36})=38$, $R(G_{37})=35$, $R(G_{38})=36$,
$R(G_{39})=33$.
\end{prop}

\begin{proof}
Note that if $n\in\{36,38\}$, then $R_G(H\setminus A)=5$. If
$n\in\{37,39\}$ then $R_G(H\setminus A)=4$. We have also
$R_G(M_1\setminus H)=4$ for $n\in\{36,37\}$ and $R_G(M_1\setminus
H)=2$ for $n\in\{38,39\}$. Now using the formulas obtained in
Lemmas \ref{fam7-A-classes}, \ref{fam7-M2} and \ref{fam7-M3-R} we
get the formulas of the proposition.
\end{proof}

\vspace{12pt} \centerline{\bf Acknowledgements}

The second author wishes to thank Mianowski Foundation (Poland)
and University of Bia{\l}ystok for the support of his research
scholarship at the University of Bia{\l}ystok, and to the first
author for his very nice hospitality during this visit. A part of
the work over the paper was done while the second author's stay in
the University of Debrecen, Hungary in March-April 2002. He would
like to express acknowledgements to Prof. Dr. Adalbert Bovdi and
Prof. Dr. Victor Bovdi for their hospitality and fruitful
discussions, and to NATO Science Fellowship Programme for the
support of this visit. We are also grateful to Prof. Dr. Bettina
Eick for her helpful advices while using the small groups library of
the computer algebra system GAP \cite{GAP4} for formulating and
testing conjectures, and to the referee for careful reading and
useful comments.

\vspace{12pt}

Czes{\l}aw Bagi\'nski \par
Institute of Computer Science \par
Technical University of Bia{\l}ystok \par
Wiejska 45A, 15-351 Bia{\l}ystok \par
Poland \par
E-mail : baginski@ii.pb.bialystok.pl

\vspace{6pt} Alexander Konovalov \par Department of Mathematics
\par Zaporozhye State University \par P.O.Box 1317, Central Post
Office, 69000 Zaporozhye, Ukraine \par E-mail: gap@gap.zssm.zp.ua,
konovalov@member.ams.org \par
http://ukrgap.exponenta.ru/konoval.htm

\end{document}